\documentclass[twocolumn]{revtex4}
\usepackage{amsmath,amsfonts,amsthm,graphics,epsfig,hyperref}

\graphicspath{{figures/}}

\newtheorem{theorem}{Theorem}
\newtheorem{lemma}[theorem]{Lemma}

\renewcommand\Pr{\operatorname{\mathbb P{}}}

\newcommand\bb[1]{\bigl(#1\bigr)}

\newcommand\Bi{\operatorname{Bi}}
\newcommand\La{\Lambda}
\newcommand\Lad{\Lambda^\star}
\newcommand\Lax{\Lambda^{\times}}
\newcommand\pc{p_{\mathrm c}}
\newcommand\pcs{p_{\mathrm c}^{\mathrm s}}
\newcommand\pcb{p_{\mathrm c}^{\mathrm b}}
\newcommand\Z{{\mathbb Z}}
\newcommand\R{{\mathbb R}}
\newcommand\cP{{\tilde \Pr}}
\newcommand\dd{{\mathrm d}}
\begin{document}

\title{Rigorous confidence intervals for critical probabilities}

\author{Oliver Riordan}
\affiliation{
 Department of Pure Mathematics and Mathematical Statistics,
 Cambridge CB3 0WB, UK.
}
\thanks{Royal Society research fellow}
\author{Mark Walters}
\affiliation{Peterhouse, Cambridge, CB2 1RD, UK}


\begin{abstract}
We use the method of Balister, Bollob\'as and Walters~\cite{BBWsquare} to give rigorous 99.9999\%
confidence intervals for the critical probabilities for site and bond percolation
on the 11 Archimedean lattices. In our computer calculations, the emphasis
is on simplicity and ease of verification, rather than obtaining the best possible
results. Nevertheless, we obtain intervals of width at most 0.0005 in all cases.
\end{abstract}

\maketitle

\section{Introduction}

In this paper we study site and bond percolation on planar
lattices, in particular the {\em Archimedean lattices},
in which all faces are regular
polygons and all vertices are equivalent. The $11$ Archimedean
lattices are shown in Figure~\ref{fig_archim}, labelled
with the
notation of Gr\"unbaum and Shephard~\cite{GS87}: each lattice
is represented by a sequence listing the numbers of sides
of the faces meeting at a vertex, in cyclic order around that vertex.

\newcommand\hh{0.7in}
\begin{figure}[htb]
\[
 \begin{array}{ccc}
  \epsfig{file=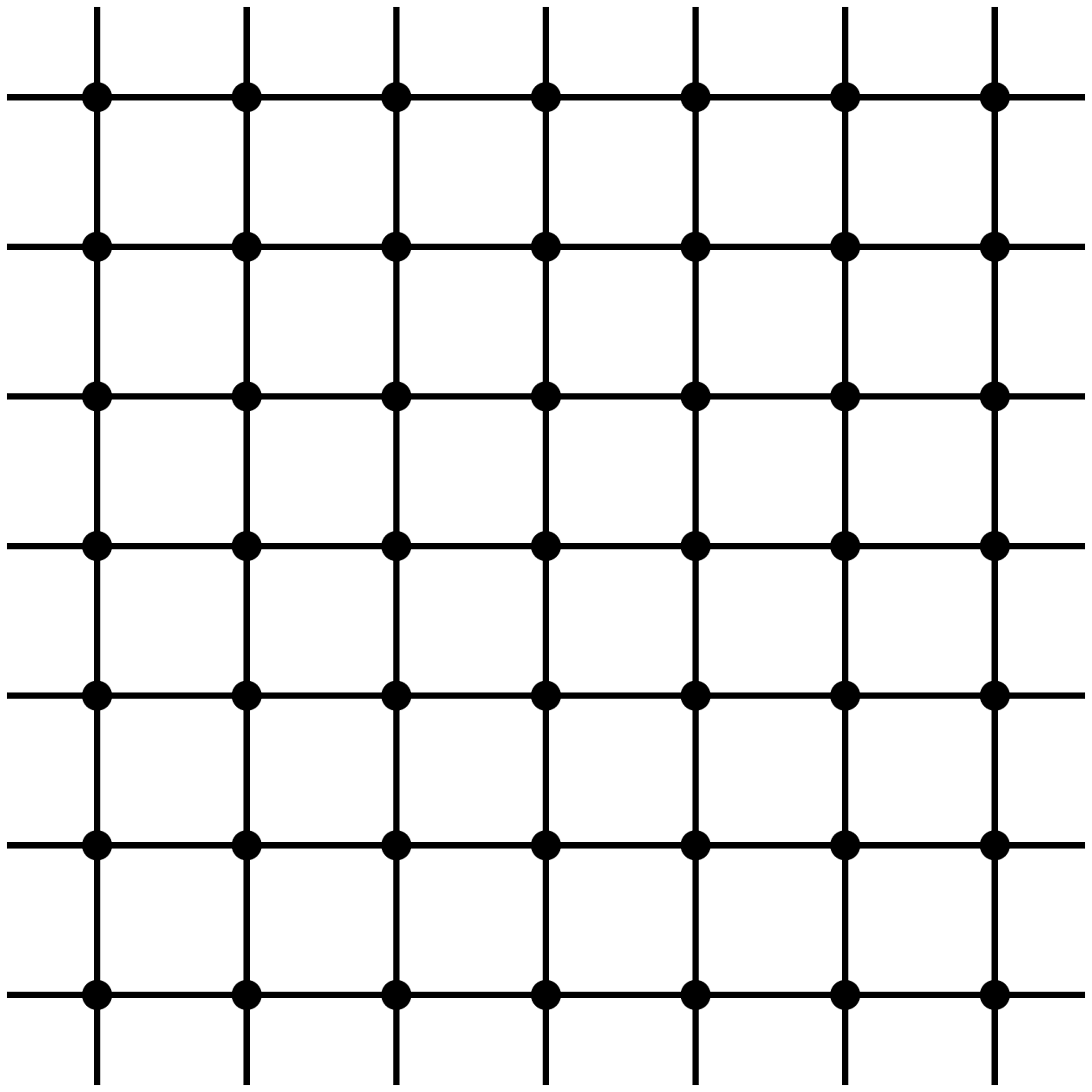,height=\hh}  &
  \epsfig{file=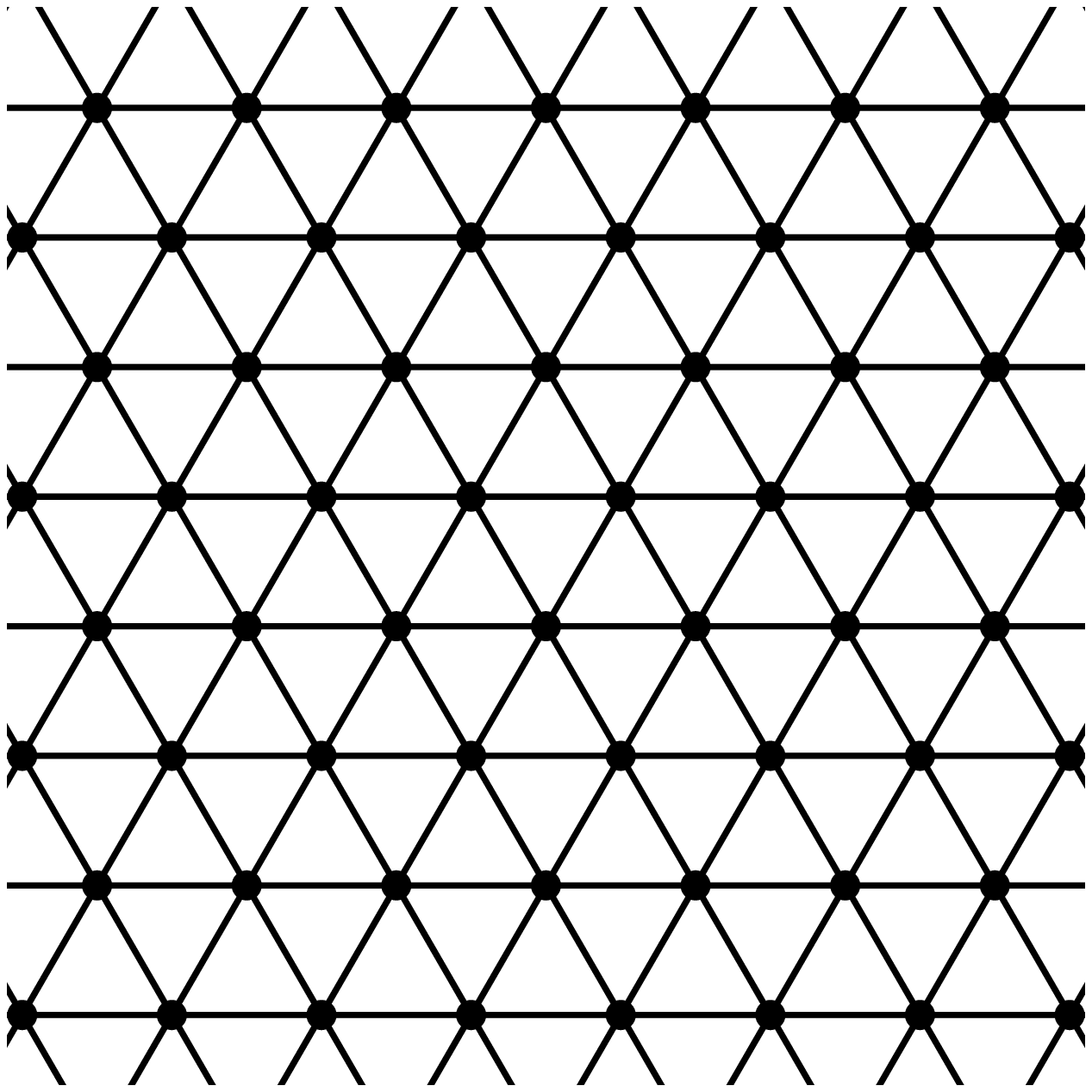,height=\hh}  &
  \epsfig{file=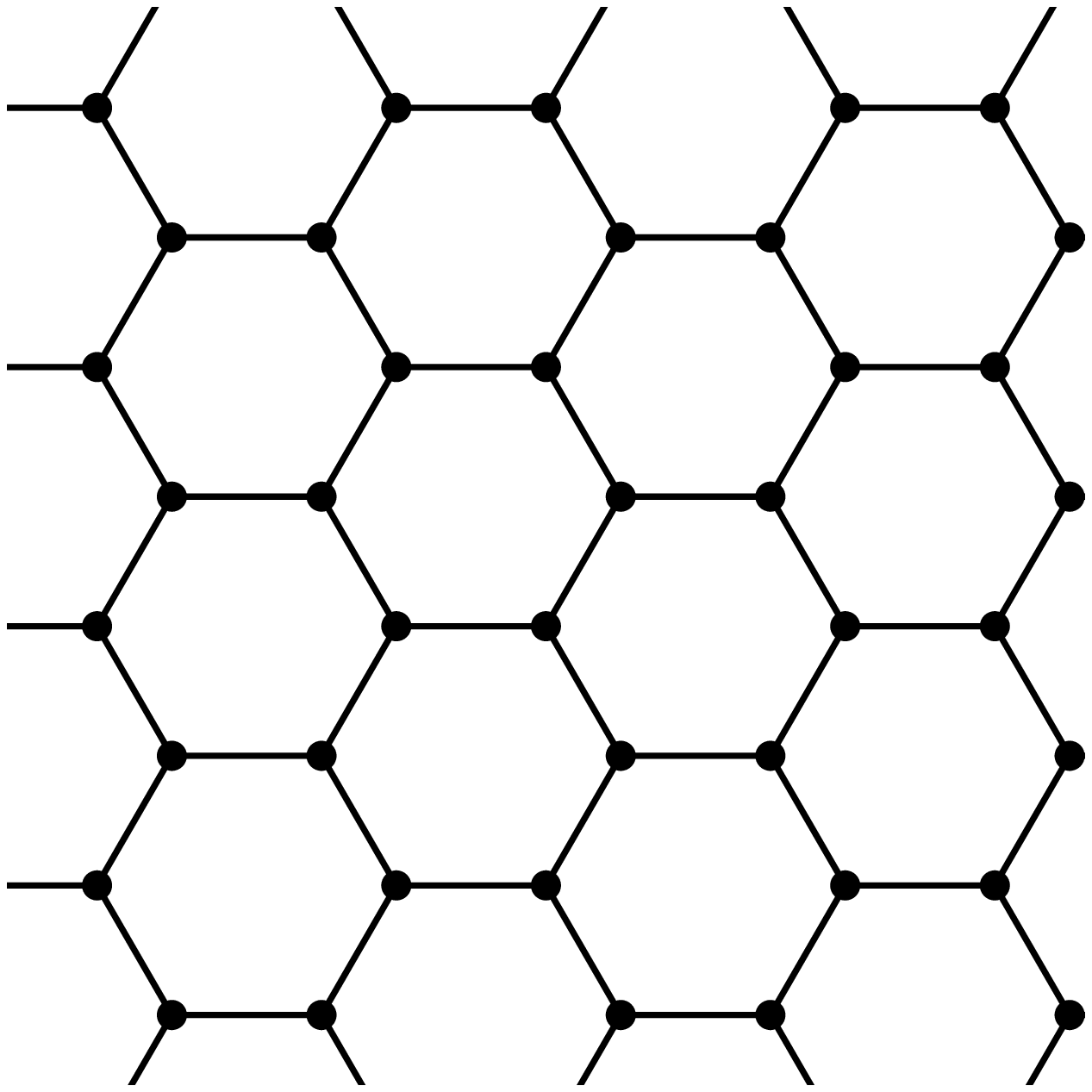,height=\hh}  \\
 \hbox{Square: }(4^4) & \hbox{Triangular: }(3^6) & \hbox{Hexagonal:
  }(6^3)  \\ \\
  \epsfig{file=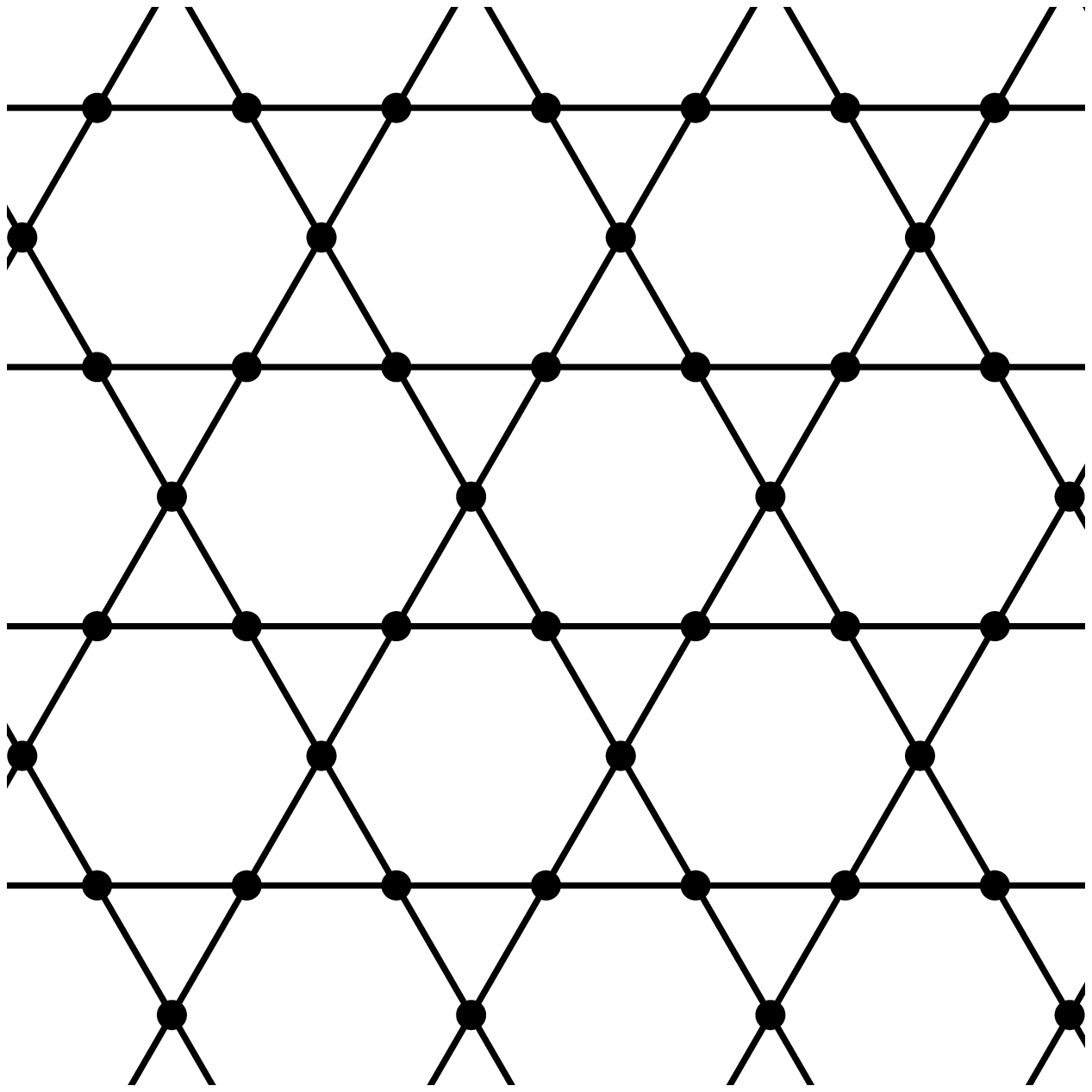,height=\hh} &
  \epsfig{file=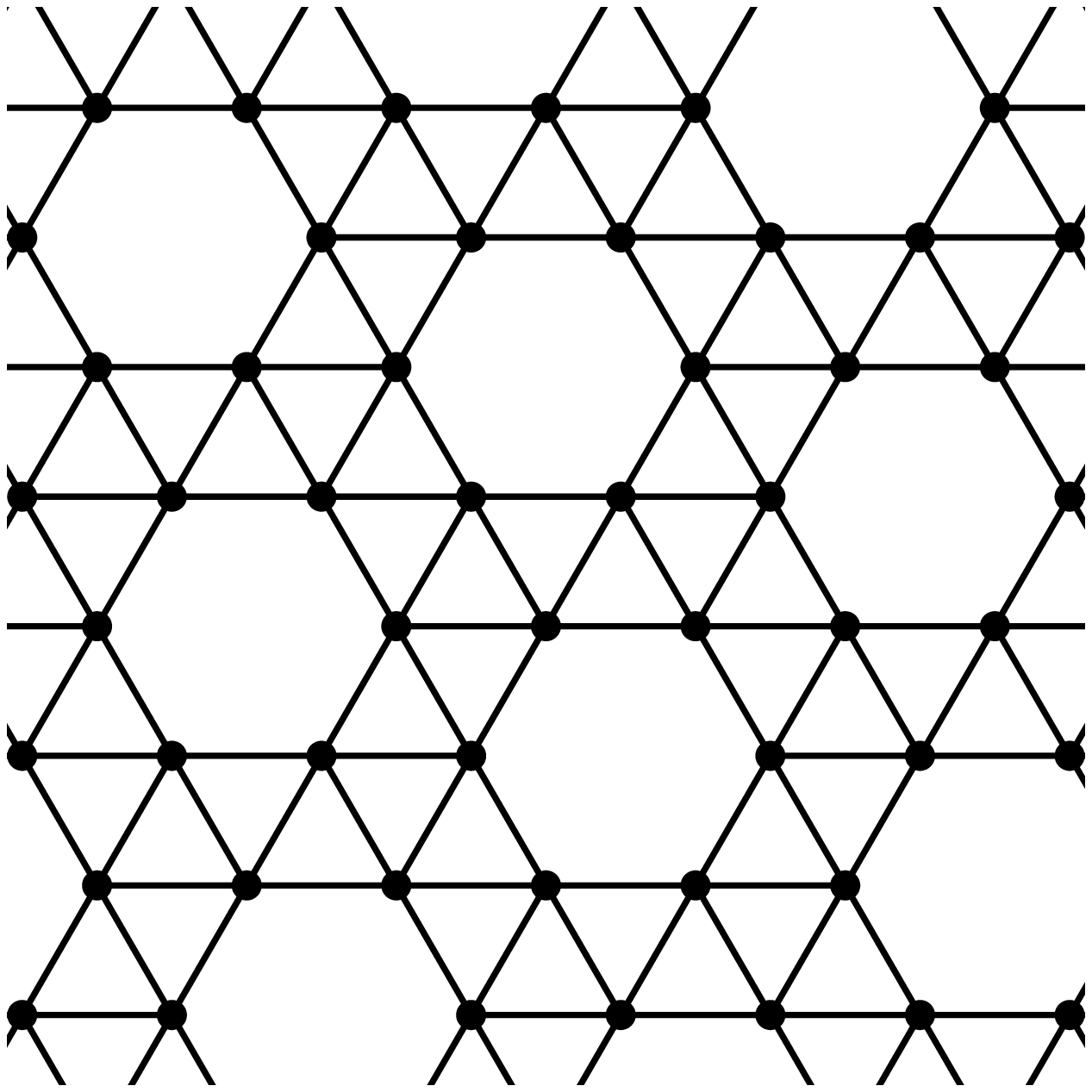,height=\hh}  &
  \epsfig{file=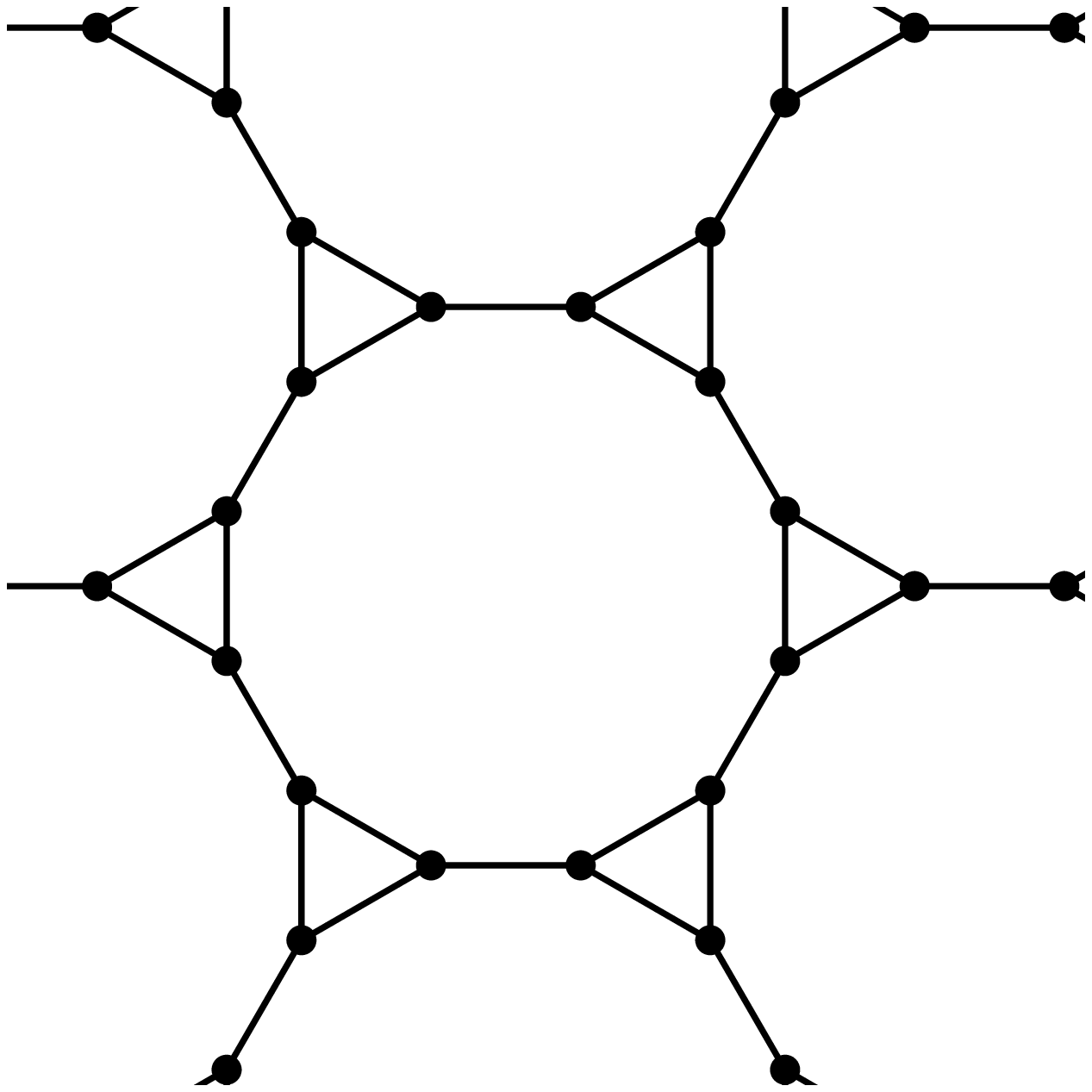,height=\hh}  \\
  \hbox{Kagom\'e: }(3,6,3,6) & (3^4,6) & (3,12^2) \\ \\
  \epsfig{file=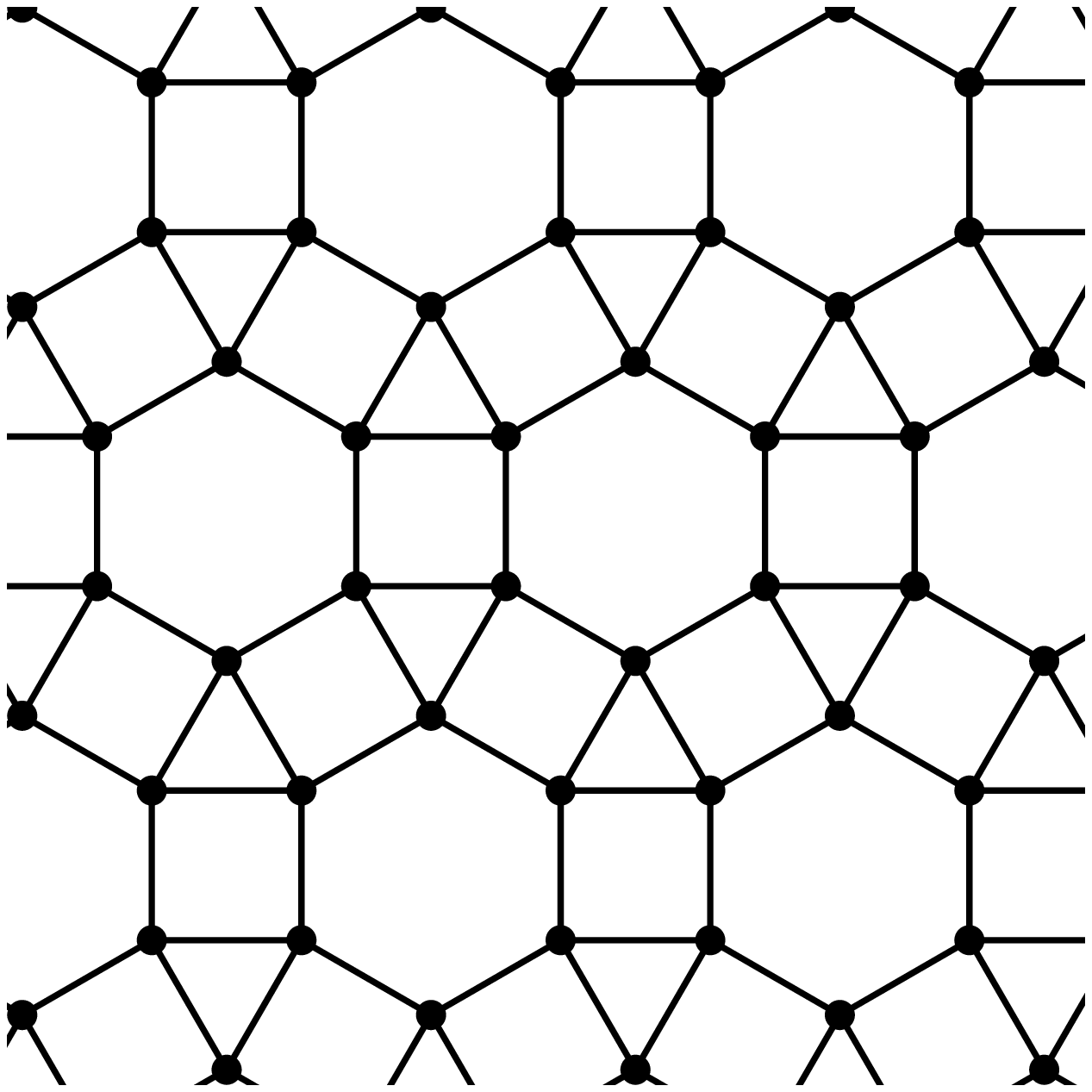,height=\hh} &
  \epsfig{file=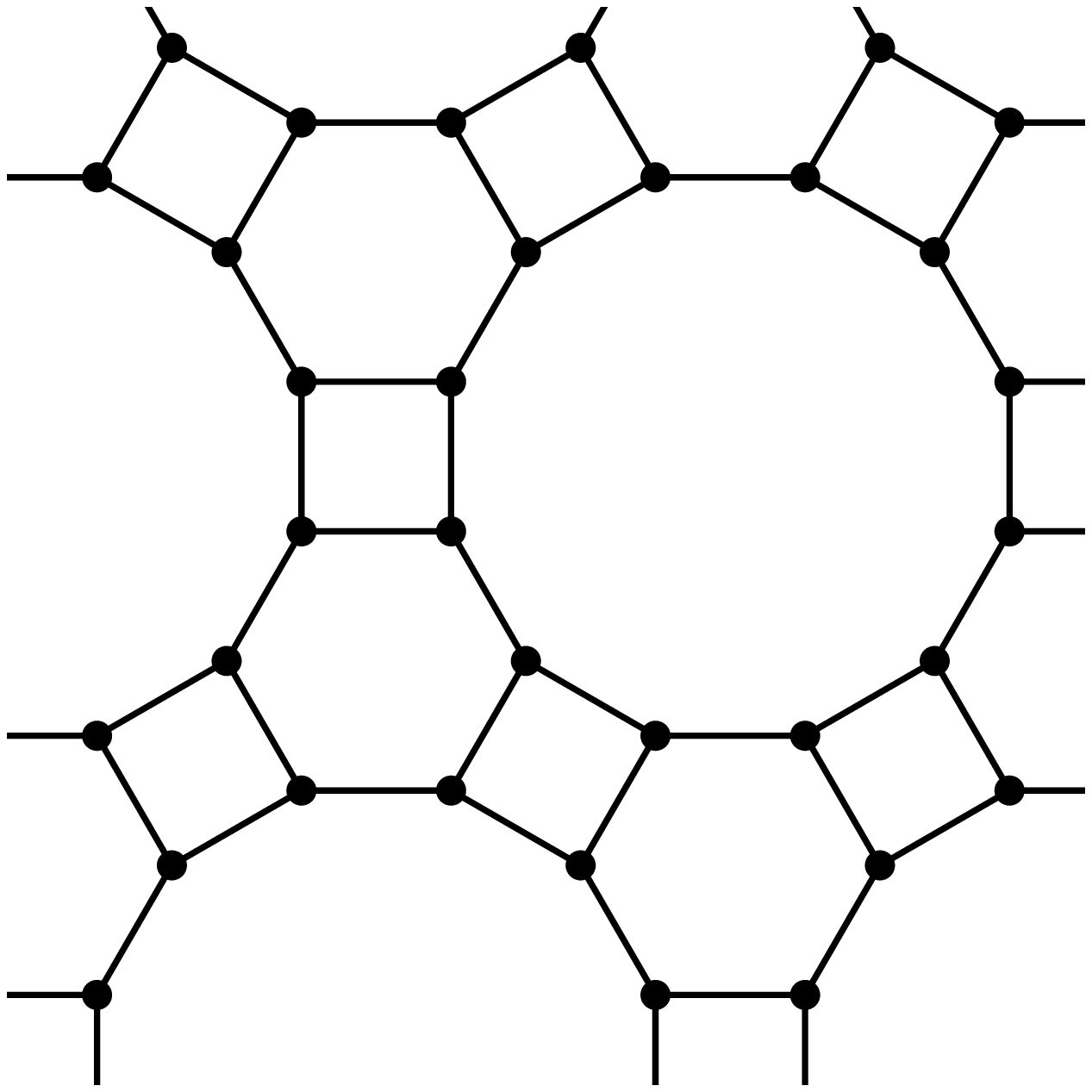,height=\hh}  &
  \epsfig{file=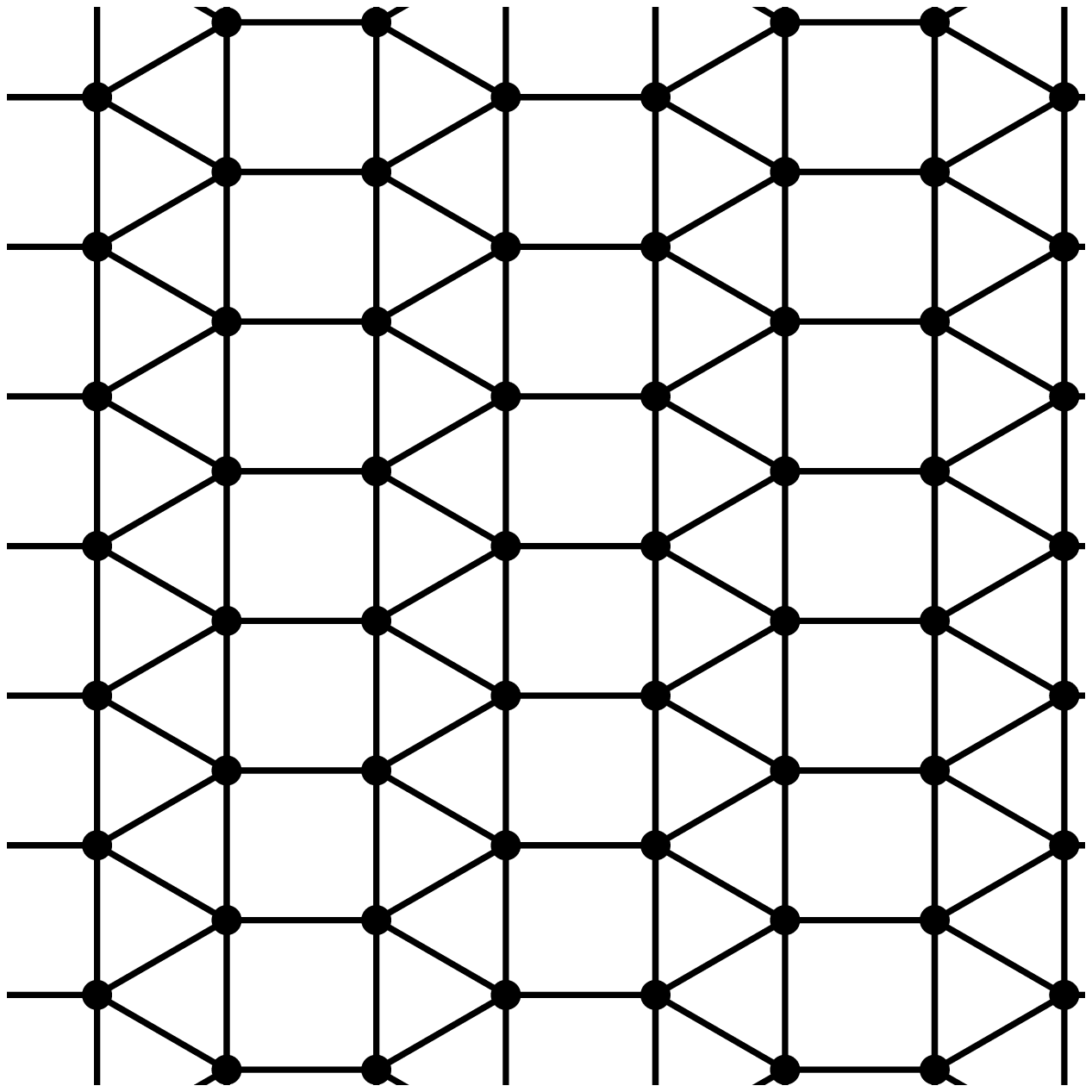,height=\hh}  \\
  (3,4,6,4) & (4,6,12)  & (3^3,4^2)\\ \\
  \epsfig{file=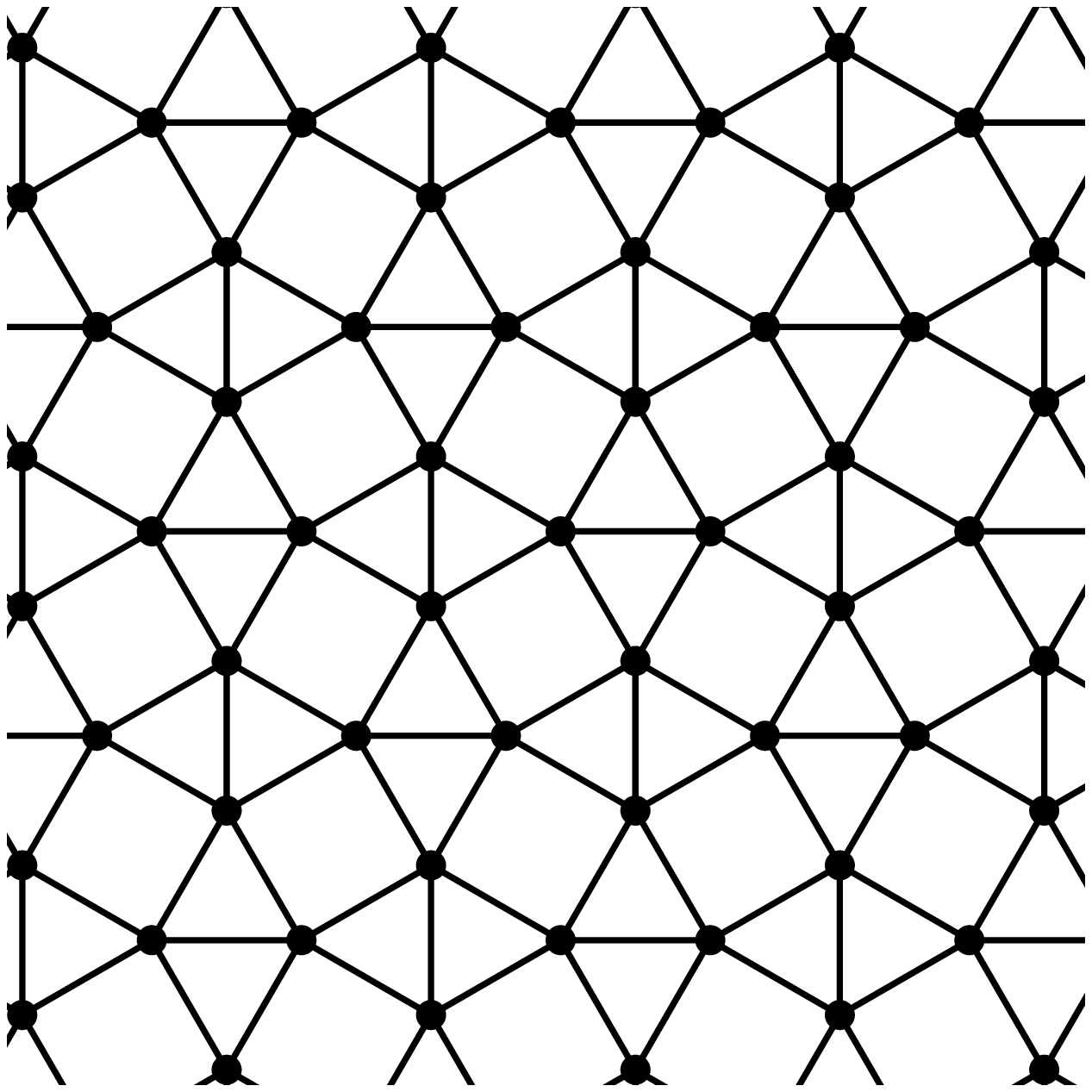,height=\hh}  &
  \epsfig{file=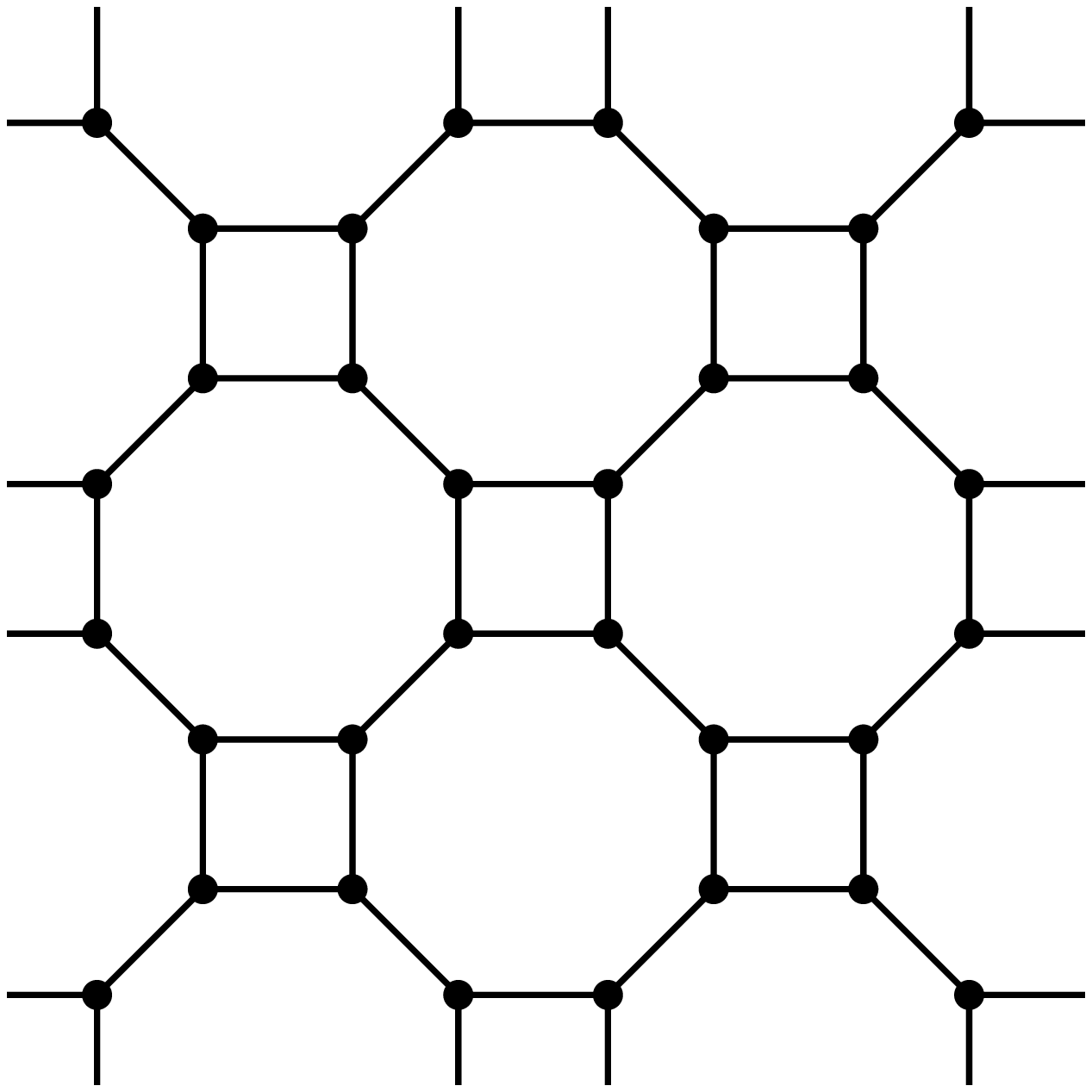,height=\hh}  \\
  (3^2,4,3,4) & (4,8^2) \\ \\
 \end{array}
\]
\caption{The 11 Archimedean lattices}
\label{fig_archim}
\end{figure}

In site percolation on a graph $\La$,
each vertex, or {\em site}, of $\La$ is assigned a {\em state}, {\em
open} or {\em closed}.  In {\em independent site percolation} the
states of the sites are independent, and each site is open with a
certain probability $p$.  The definitions for bond percolation are
similar, except that is the edges, or {\em bonds}, of $\La$ that are
assigned states. We shall write $\Pr_p$ for the corresponding
probability measure, suppressing the dependence on $\La$ and on
whether it is site or bond percolation that we consider.

The basic question of percolation theory is `when is there an infinite
{\em open cluster}', i.e., an infinite subgraph of $\La$ all of whose
sites (for site percolation) or bonds (for bond percolation) are open.
It is not hard to see that there is a certain `critical probability'
$\pc$, such that for $p<\pc$ there is never (i.e., with probability
$0$) an infinite open cluster, while for $p>\pc$ there always is.  For
this and other basic facts about percolation, see
Grimmett~\cite{Grimmett2nd}, or Bollob\'as and Riordan~\cite{BRbook}, for example.
When we wish to
specify the lattice, and whether it is site or bond percolation that
we are considering, then we write $\pcs(\La)$ or $\pcb(\La)$.

The exact value of $\pc$ is known in rather few cases:
in 1980, Kesten~\cite{Kesten1/2} proved that $\pcb(\Z^2)=1/2$,
where $\Z^2$ is the square lattice.
Shortly afterwards~\cite{KestenBook}, he proved that $\pc=1/2$ also holds for site
percolation on the triangular lattice $T$. Later, Wierman~\cite{WiermanHoneycomb}
used his `substitution' method to give rigorous proofs of the values
$\pcb(T)=2\sin(\pi/18)$ and $\pcb(H)=1-2\sin(\pi/18)$ for bond percolation on the
triangular and hexagonal lattices respectively; these values had been
obtained heuristically much earlier by Sykes and Essam~\cite{SykesEssamPRL,SykesEssam64}.
There are two further values that may be easily derived from these:
the $(3,6,3,6)$ or {\em Kagom\'e} lattice $K$ is the line graph of the hexagonal lattice, so
$\pcs(K)=\pcb(H)=1-2\sin(\pi/18)$. Also,
the $(3,12^2)$ or {\em extended Kagom\'e} lattice $K^+$ is the line graph
of the lattice $H_2$ obtained by subdividing each bond of $H$ exactly once,
so
\[
 \pcs(K^+)=\pcb(H_2)=\sqrt{\pcb(H)}=\bb{1-2\sin(\pi/18)}^{1/2}.
\]
These are the only critical probabilities known for Archimedean lattices.
Indeed, it may be
that the exact values of the other critical probabilities associated
to the Archimedean lattices will never be known; they may simply be numbers
that have no simpler descriptions than their definitions as critical probabilities.

Given the dearth of exact results, it is not surprising that much effort has
been put into the estimation of critical probabilities. Almost all results
in this area are of one of two types: (1) rigorous upper and/or lower bounds,
and (2) heuristic estimates based on computer calculations. There are also a few
heuristic derivations of conjectured exact results; we shall
return to this briefly in Section~\ref{concl}. Examples
of (1) are the bounds obtained by Wierman~\cite{Wierman_sub1,Wierman_Kagome2003,WiermanArchim}
(see also Parviainen and Wierman~\cite{PW}) using his substitution
method. Even with considerable work on efficient algorithms and extensive computer
calculations, it seems to be hard to obtain narrow intervals between rigorous
upper and lower bounds: of the intervals listed in~\cite{PW}, two have
width a little under $0.01$, but many have width $0.1$ or even $0.2$.

There are a small number of recent, more accurate, rigorous results,
including an intervals of width $0.00135$ and $0.0046$ for site percolation on
the $(3,12^2)$ and Kagom\'e lattices obtained by May and Wierman \cite{MayWierman}.

Turning to (2), there are so many papers on this topic, going back to the 1960s,
that it is impossible to attempt even a representative list. 
Let us mention a couple of
examples, however: for $\pcs(\Z^2)$, Reynolds, Stanley and Klein \cite{Reynolds}
reported the impressively accurate non-rigorous estimate $0.5931\pm
0.0006$ already in 1980. In 1986, Ziff and Sapoval \cite{ZiffSapoval} gave the
exceedingly accurate estimate of $0.592745(2)$. In 1990 Yonezawa,
Sakamoto and Hori \cite{Yonezawa} gave estimates for several
Archimedean lattices, with error terms a little over $10^{-4}$.
More recently, many very precise  estimates have been given
using quite sophisticated methods, for example by Suding and Ziff~\cite{SudingZiff},
Newman and Ziff~\cite{NewmanZiff_fast,NewmanZiff_efficient}  
and Parviainen~\cite{Parv_est}. It is very likely that these
estimates are extremely accurate; errors of `about $\pm 3\times 10^{-6}$'
are claimed in~\cite{SudingZiff}, and even smaller errors
in~\cite{NewmanZiff_fast,NewmanZiff_efficient,Parv_est}.
However, these estimates come with no mathematical guarantees, and it is hard to be
sure how accurate they really are. Although theoretical
error analysis is sometimes given (see Ziff and Newman~\cite{ZiffNewman_scaling}, for
example), this is certainly non-rigorous. Even assuming unproved
results about the scaling limits of planar percolation models gives only the asymptotic
behaviour of these errors; it does not allow us to say anything about the relationship
between a finite number of data points and the true value of $\pc$. Also, there is
disagreement about even the asymptotic form of the errors in some cases (see
Parviainen~\cite{Parv_est}), and there are several
instances where earlier estimates have been contradicted by later ones.

Surprisingly, it is possible to give a result intermediate in nature between
a rigorous bound and a heuristic estimate: one can {\em prove} that a certain (random)
procedure generates a bound that is correct with probability at least 99.9999\%, say;
in other words, one can rigorously generate confidence intervals for critical probabilities.
Such intervals are typically much narrower than the 100\% bounds, although nothing
like as narrow as the (claimed) uncertainties for heuristic estimates.
Results of this kind were first proved by Bollob\'as and Stacey~\cite{BS_oriented} in the context
of oriented percolation, and then by Balister, Bollob\'as and Walters~\cite{BBWsquare}
in the context of (unoriented) continuum percolation.
Here we use the method of the latter paper, which applies essentially
`as is' to percolation on 2-dimensional lattices, to obtain confidence intervals
for the site and bond percolation critical probabilities for all 11 Archimedean
lattices. Indeed such an interval for site percolation on the square lattice
was given in~\cite{BBWsquare}, namely $[0.5919,0.5935]$ (with a lower
confidence of 99.99\%). Here, with greater 
computational effort, we obtain narrower intervals.

\section{Method}

\subsection{The mathematics}\label{subsec_math}

The method of Balister, Bollob\'as and Walters~\cite{BBWsquare} is based on
a simple application of the concept of {\em 1-independent percolation}
(also known as {\em 1-dependent percolation}).  A bond percolation
measure on a graph $\La$, i.e., a measure on assignments of states to the bonds
of $\La$, is {\em 1-independent} if, whenever $S$ and $T$ are sets of
bonds such that the graph distance from $S$ to $T$ is at least $1$,
the states of the bonds in $S$ are independent from the states of
the bonds in $T$. In other words, roughly speaking,
the states of vertex-disjoint bonds are independent.
Such measures arise naturally in percolation theory, in particular in
static renormalization arguments (see Section 7.4 of Grimmett~\cite{Grimmett2nd}),
and have been considered by many authors.
Although the assumption of $1$-independence is weaker than independence, it is strong enough
to ensure percolation if the individual bonds are open with sufficiently high probability,
as shown by the following lemma of Balister, Bollob\'as and Walters~\cite{BBWsquare}.

\begin{lemma}\label{l_BBW}
Let $\cP$ be a $1$-independent bond percolation measure on $\Z^2$ in which
each bond is open with probability at least $p_0=0.8639$. Then the
probability that the origin lies in an infinite open cluster is positive.
\end{lemma}

If the value of $p_0$ is not important, then a weak form of Lemma~\ref{l_BBW}
(with $p_0$ replaced by some constant smaller than $1$) is more or less
immediate from first principles (see Bollob\'as and Riordan~\cite{BRKesten}, for example). It also
follows from the very general results of Liggett, Schonmann and Stacey~\cite{LSS}
comparing $1$- (or $k$-) independent measures on general graphs with product measures.

Starting from independent site or bond percolation on a lattice
$\La\subset \R^2$, there is a natural way to obtain a $1$-independent
bond percolation measure $\cP$ on the square lattice $\Z^2$:
given a `scale parameter' $s>0$, partition $\R^2$ into
disjoint $s$ by $s$ squares $S_v$, $v\in \Z^2$. For each bond $e$ of
$\Z^2$ let $R_e$ be the corresponding rectangle, so if $e=uv$ then
$R_e=S_u\cup S_v$.  Let $E_e$ be some event that depends only on the
states of the sites or bonds of $\La$ that lie within $R_e$, and take
the bond $e$ of $\Z^2$ to be open with respect to $\cP$ if and only if $E_e$ holds.  Since the
rectangles corresponding to vertex-disjoint bonds of $\Z^2$ are
disjoint, this defines a $1$-independent measure~$\cP$.

Suppose that the squares $S_v$ and events $E_e$ are chosen so that
the following condition holds:
\begin{quotation}
 whenever there is an infinite path $v_0v_1v_2\ldots$
 such that $E_{e_i}$ holds for each $e_i=v_iv_{i+1}$, there
is an infinite open cluster in the original lattice
\end{quotation}
\vskip-0.62in \begin{equation}\label{Acond}\end{equation}\vskip0.1in
\noindent
and then $p$ is chosen so that
\begin{equation}\label{pcond}
 \Pr_p(E_e)\ge 0.8639 \hbox{ for every bond $e$ of $\Z^2$}.
\end{equation}

Then Lemma~\ref{l_BBW} implies that $\pc\le p$. Indeed, we have already
noted that $\cP$ is $1$-independent, so from \eqref{pcond} and Lemma~\ref{l_BBW}
there is a positive probability that the origin is in an infinite $\cP$-open
cluster. But then (as $\Z^2$ is locally finite) there is an infinite $\cP$-open path
starting at the origin, so from \eqref{Acond} there is an infinite open cluster
in the original percolation with positive probability. Hence $p\ge \pc$, as required.
This type of argument, but in a qualitative form where the value
of $p_0$ is not important, provides one of the many easy ways of deducing
Kesten's $\pcb(\Z^2)=1/2$ result from a suitable `sharp-threshold' result;
see Bollob\'as and Riordan~\cite{BRKesten,BRbook}.

Here we follow Balister, Bollob\'as and Walters~\cite{BBWsquare} in our choice 
for the event $E_e$. For $v\in \Z^2$, let $\La_v$ denote the subgraph of $\La$
induced by the sites in $S_v$. For each bond $e=uv$ of $\Z^2$, let
$\La_e$ by the subgraph of $\La$ induced by the sites in $R_e=S_u\cup S_v$.
Let $E_e$ be the event that each of $\La_u$ and $\La_v$ contains a unique largest open cluster
with these clusters part of the same open cluster in $\La_e$.
Here `largest' simply means containing the most
sites. Note that $E_e$ does depend only on the states of bonds or sites
within $R_e$, so we do obtain a $1$-independent measure. Also,
it is immediate that \eqref{Acond} is satisfied.

To obtain an upper bound on $\pc$, it remains only to find a pair $(s,p)$ for which \eqref{pcond} is satisfied.
Note that it will suffice to check condition \eqref{pcond} for (usually) one
or (occasionally) two bonds $e$ of $\Z^2$: without changing the
graph structure, we shall redraw all our lattices $\La$ so that the vertex
set is a subset of $\Z^2$, and so that horizontal and vertical translations
through some small integer $C$ act as isomorphisms of $\La$.
For $u=(a,b)$, $a,b\in \Z$, we shall take
\[
 S_u = \{(x,y): sa\le x<s(a+1),\,sb\le y<s(b+1)\},
\]
where $s$ is a `scale' parameter with $C$ dividing $s$. Thus all squares $S_v$
are equivalent with respect to the lattice. Furthermore, for the lattices
with an axis of symmetry, our new representation will have the line $x=y$
as an axis of symmetry. This ensures that all rectangles $R_e$ are equivalent,
so $\Pr_p(E_e)=\Pr_p(E_f)$ for all $e$, $f$. When there is no axis of symmetry,
we have to consider one rectangle with each orientation.

So far, we have only discussed upper bounds; this is because we can
obtain lower bounds by bounding the critical probability for a related
lattice from above. Indeed, given a planar lattice $\La$, let $\Lad$
be the usual planar dual of $\La$, with one site for each face of
$\La$, and a bond $e^\star$ for each bond $e$, joining the two sites of
$\Lad$ corresponding to the faces in which $e$ lies. It is `well known'
that
\begin{equation}\label{dsum}
 \pcb(\La) + \pcb(\Lad) =1.
\end{equation}
Thus, to bound $\pcb(\La)$ from below we may bound $\pcb(\Lad)$ from above.

Note that while \eqref{dsum} is widely assumed to be true in great
generality, it has only been proved under certain symmetry assumptions.
Under very general conditions, the upper bound $\pcb(\La)+\pcb(\Lad)\le 1$
follows immediately from Menshikov's Theorem~\cite{Menshikov}.
For the lower bound, one shows that it is not possible to have bond percolation
in $\La$ at a parameter $p$ and also bond percolation in $\Lad$ at parameter $1-p$.
For lattices (doubly periodic, locally finite planar graphs) with
rotational symmetry of some order $k\ge 4$, there is a simple proof of the lower bound
due to Zhang; see Lemma 11.12 of Grimmett~\cite{Grimmett2nd}, where this argument is presented
for $\Z^2$. Recently, Bollob\'as and Riordan~\cite{BRdual} (see also~\cite{BRbook})
have pointed out that
this argument can be easily adapted to lattices with rotational symmetry of any order $k\ge 2$.
This is important here: all 11 Archimedean lattices have such rotational symmetry,
but two, the lattices $(3^2,4,3,4)$ and $(3^3,4^2)$, do not have rotational symmetry of
higher order. Even more recently, Sheffield~\cite{Sheffield} has given a much more complicated
argument that proves \eqref{dsum} for lattices without further symmetry assumptions.

For site percolation, let $\Lax$ be the (in general non-planar) graph obtained from $\La$
by adding a bond between any two sites in the same face of $\La$; we shall
refer to $\Lax$ as the {\em site dual} of $\La$. One has
\[
 \pcs(\La) + \pcs(\Lax)=1;
\]
the comments above about symmetry assumptions apply in this case also.

\subsection{The statistics}
For sufficiently small scale parameters $s$, it is possible to find a
$p$ for which \eqref{pcond} holds by enumerating all possibilities
for which sites/bonds in $R_e$ are open, and so writing $\Pr_p(R_e)$ as
a polynomial in $p$. Needless to say, this is impractical and gives poor results
in practice.
The key idea of Balister, Bollob\'as and Walters~\cite{BBWsquare} is to use a statistical
approach, obtaining confidence intervals with precisely calculated error probabilities
instead of 100\% upper bounds.
Indeed, suppose that we have a random procedure $A$ for generating a pair $(s_A,p_A)$,
and that one can {\em prove} that, with probability at least $99.9999\%$,
the random pair produced is one for which \eqref{pcond} holds. Then
$(-\infty,p_A]$ is a (random, as always) one-sided $99.9999\%$ confidence interval
for $\pc$ (see below). Such a procedure $A$ is very easy to define;
again, we follow \cite{BBWsquare}, with one small modification (and with
different numbers).

Suppose that we have somehow `guessed' values of the scale parameter $s$
and percolation parameter $p$ for which we expect
that $\Pr_p(E_e)$ is somewhat larger than $0.8639$. We then generate
$N=400$ random simulations of the configuration within $R_e$,
and count the number $m$ of them in which $E_e$ holds.
If $\Pr_p(E_e)=\pi$, then $m$ has a binomial $\Bi(N,\pi)$ distribution with parameters
$N$ and $\pi$.
In particular, {\em if} $\pi<0.8639$, then
\begin{align*}
 \Pr(m\ge 378) &\le \Pr\bb{\Bi(400,0.8639)\ge 378} \\
&= 1.1489\dots\times 10^{-7} < 10^{-6}/6.
\end{align*}
If our simulation does give $m\ge 378$, we can thus assert with very high confidence
that $\pi\ge 0.8639$, i.e., that \eqref{pcond} does hold, which, as noted
above, implies $\pc\le p$.

This is the heart of the method of Balister, Bollob\'as and Walters~\cite{BBWsquare} (and 
also of the related method of Bollob\'as and Stacey~\cite{BS_oriented}): no matter how we arrive
at our `guess' for $s$ and $p$, provided we only perform one `final' simulation, the simple
inequality above shows that the probability that we assert an incorrect upper bound
for $\pc$ is at most $10^{-6}/6$. Note that we may be unlucky: if $m<378$,
then we can assert only the trivial upper bound $1$. In terms of the description above,
our random procedure returns the guessed values $(s,p)$ if $m\ge 378$, and the trivial
pair $(s,1)$ otherwise. Note that we have no bound on the probability that we get
$1$ as an upper bound, but this does not matter for the argument that $(-\infty,p_A]$
is a $99.9999\%$ confidence interval. Of course, to obtain useful results,
we want to be reasonably sure that we will have $m\ge 378$, and this is where
the careful choice of parameters comes in.

Here we modify the method very slightly: the choice of the number $378$ gives us individual
error probabilities that are smaller than $10^{-6}/6$. Hence, we can perform up to three
different runs with different parameters $s$ and $p$ (which may depend on the results of previous
runs), and choose the best bound
given by a successful run. It is still true that each run has at most a probability $10^{-6}/6$
of producing an incorrect bound, so the probability that our final bound is incorrect is at most
$10^{-6}/2$. Bearing in mind that the same applies to the lower bounds (realised as
upper bounds on a dual critical probability), we still obtain $99.9999\%$ confidence intervals.

A small side note: since the lattice $(3^4,6)$ does not have an axis of
symmetry, we ran our method in two directions, horizontally
and vertically. This means that we wanted the individual error probability 
to be smaller than $10^{-6}/12$ which was satisfied
by requiring at least $379$ successes in this case.

There are two advantages to this method: it turns out to be slightly
more efficient (based on heuristic calculations). Bearing in mind that
we can stop after one successful run, we can perform three runs each
of which has a 90\% chance of succeeding, say, more quickly than one
run for the same $p$ but a larger $s$ that has $99.9\%$ chance of
succeeding. Secondly, if we are not very confident of our guesses,
after a failed first run we can choose more conservative parameters for the second
and third runs (for example, keeping $p$ fixed but increasing $s$),
to be very sure of obtaining reasonable bounds in the end.

\subsection{Random number generation}

So far, we have assumed the availability of a suitable source of
random numbers. In practice, one usually uses a pseudo-random number
generator. This introduces a possible source of error: it could be
that there is some pattern in the output of the generator that affects
the results of the simulations.  To minimize the likelihood of this we
used the well known and well trusted MT19937 ``Mersenne Twister''
generator developed by Matsumoto and Nishimura \cite{MT1}, as
updated in 2002. See their website \cite{MTweb} for the source code
and related literature.

It would be very easy to modify our program to use other random
number generators, or even a hardware generator.

The selection of a random number generator for simulations is often glossed
over; here we emphasize this as it is important for our results:
the {\em only} assumption in our results (that our procedure produces
99.9999\% confidence intervals for $\pc$) is that the
random numbers used in the simulation may be treated as genuinely random. 

\subsection{Choice of parameters}

In this subsection we outline the purely heuristic arguments we used to choose
suitable parameters for running our final statistical tests. The correctness
of the results does not depend on the correctness of these arguments. For
this reason we allow ourselves to use consequences of the very widely believed but, except
for one lattice, unproved conformal invariance conjecture.
This conjecture of Aizenman and Langlands, Pouliot and Saint-Aubin~\cite{Langlands_confinvar}
states (among other things) that, for any planar lattice,
after a suitable affine transformation, the limiting crossing probabilities
for large regions are invariant under conformal mappings, and, more precisely,
are given by Cardy's formula~\cite{Cardy92}. For more details see Bollob\'as and Riordan~\cite{BRbook}, for example.
As shown by Smirnov and Werner~\cite{SmirnovWerner}, building on work of
Schramm~\cite{SchrammSLE} and 
Lawler, Schramm and Werner~\cite{LSW_bie1,LSW_bie2,LSW_onearm,LSW_bie3},
this conjecture, if true, enables the values of certain `critical exponents' to be calculated.
Note that the conjecture has been proved, by Smirnov~\cite{Smirnov}, only for site percolation
on the triangular lattice; for all other lattices it is still open.

Fixing the percolation model under consideration, i.e., fixing the lattice
$\La$, and considering either bond or site percolation throughout, let $\pc$
be the appropriate critical probability, and set
\[
 f(s,p) = \Pr_p(E_e)
\]
for one fixed bond $e$ of $\Z^2$, noting that the definition of the squares
$S_u$ and hence of the event $E_e$ depends on our scale parameter $s$.
It is not hard to convince oneself that $f(s,\pc)$ tends to some constant $0<a<1$
as $s\to\infty$, although this does not obviously follow formally from the conformal invariance
conjecture.

Turning to the $p$-dependence of $f(s,p)$, it is natural to guess that for fixed $s$,
for $p$ not too far from $\pc$, the function $f(s,p)$ will roughly satisfy the differential
equation
\[
 \frac{\dd}{\dd p} f(s,p) = C(s) f(s,p)(1-f(s,p)),
\]
where $C(s)$ is a constant depending on $s$ (and on the lattice). For one thing, $f(s,p)$ should
decay exponentially, and approach $1$ exponentially, as $p$ moves away from $\pc$.
Also, by the Margulis-Russo formula, $\frac{\dd}{\dd p} f(s,p)$ is exactly the expected
number of sites/bonds that are {\em pivotal} for the event $E_e$, i.e., such that changing the
state of this site/bond from closed to open or vice versa alters whether $E_e$ holds.
If a site (say) $v$ is pivotal, then $E_e$ must hold in the configuration with $v$ open,
and not hold with $v$ closed, so it is reasonable to guess that
for fixed $s$ and a fixed site $v$, the probability that
$v$ is pivotal will be roughly proportional to $\Pr_p(E_e) (1-\Pr_p(E_e))$.

Up to a constant factor, $C(s)$ above is just $s^2$ times the probability that
a `typical' site (or bond) $v$ is pivotal for $E_e$ at $p=\pc$. Roughly speaking,
$v$ is pivotal if and only if, when $v$ is open, two open clusters of (linear) scale $s$
are joined which, if $v$ is closed, are separated by a path of linear scale $s$.
Hence, the probability that $v$ is pivotal should scale as $s^{-\alpha_4}$, where
$\alpha$ is the `multi-chromatic $4$-arm
exponent'. Roughly speaking, $\alpha_4$ is defined as the scaling exponent
of the probability that there are four disjoint paths $P_1,P_2,P_3,P_4$ from $v$ (or from `near' $v$)
to points at distance $s$ from $v$, with $P_1$ and $P_3$ open, $P_2$ and $P_4$ closed,
with the endpoints of the $P_i$ appearing in cyclic order.
Assuming conformal invariance, from \cite{SmirnovWerner} we have $\alpha_4=(4^2-1)/12=5/4$,
so we expect $C(s)$ to scale as $s^{2-5/4}=s^{3/4}$.

Putting the above together, it is reasonable to expect the function $f(s,p)$ to have
approximately the form
\begin{equation}\label{fform}
 \frac{1}{1+\exp\bb{a-bs^{3/4}(p-\pc)}},
\end{equation}
for some constants $a$ and $b>0$ that depend on the lattice.
Our procedure for choosing the final parameters $(s,p)$ to use is as follows:
first, numerically estimate $f(s,p)$ for a fixed small value $s_0$ of $s$ (typically 72)
and various values of $p$. Then fit the data with the function above to give a rough
estimate of $a$ and $b$. Then calculate values $p_{1/3}$, $p_{2/3}$ of $p$ at which the formula
predicts $f(s_0,p_{1/3})=1/3$ and $f(s_0,p_{2/3})=2/3$.
Next, run more extensive simulations to estimate $f(s_0,p_{i/3})$, and use these
two datapoints to calculate better estimates of $a$ and $b$.
The reason for this step is that we do not expect \eqref{fform} to give a very
accurate description of
the shape of the curve $f(s,p)$ with $s$ fixed and $p$ varying, particularly
when $p$ is far from $\pc$,
so we wish to extrapolate from consistently chosen points on this curve.

Finally, we aim to choose a (large) $s$ and a $p$ close to $\pc$
such that $f(s,p)$ is approximately 0.957; this is because $\Pr\bb{\Bi(400,0.957)\ge 378}$ is close
to $90\%$, so with these parameters we have a good enough chance of obtaining a valid bound, bearing
in mind that we can perform three separate runs. Extrapolating \eqref{fform} this far does
not give very good results; experimentally, when \eqref{fform} is
about $0.945$, or a little less, the true value of $f(s,p)$ is large enough.
Of course, the larger $s$ is, the closer $p$ can be taken to $\pc$.
The exact values of $s$ and $p$ were chosen based on the amount of computer time available,
and so that we obtained intervals of width at most $0.0005$ in all cases.

\section{Computations and results}

Although our final results are confidence intervals, we are aiming for
{\em rigorous} confidence intervals, i.e., we must {\em prove} that,
for each lattice, our procedure has probability at least $99.9999\%$
of producing an interval containing the true value (assuming the random
number generator we used is well behaved).
The main
practical consequence of this is that we must ensure that we evaluate
$\Pr_p(E_e)$ for rectangles $R_e$ that fit together exactly in the
manner required for the argument in Subsection~\ref{subsec_math}.

\begin{figure}
\[
\rotatebox{-90}{\reflectbox{\rotatebox{90}{\resizebox{!}{100pt}{\includegraphics{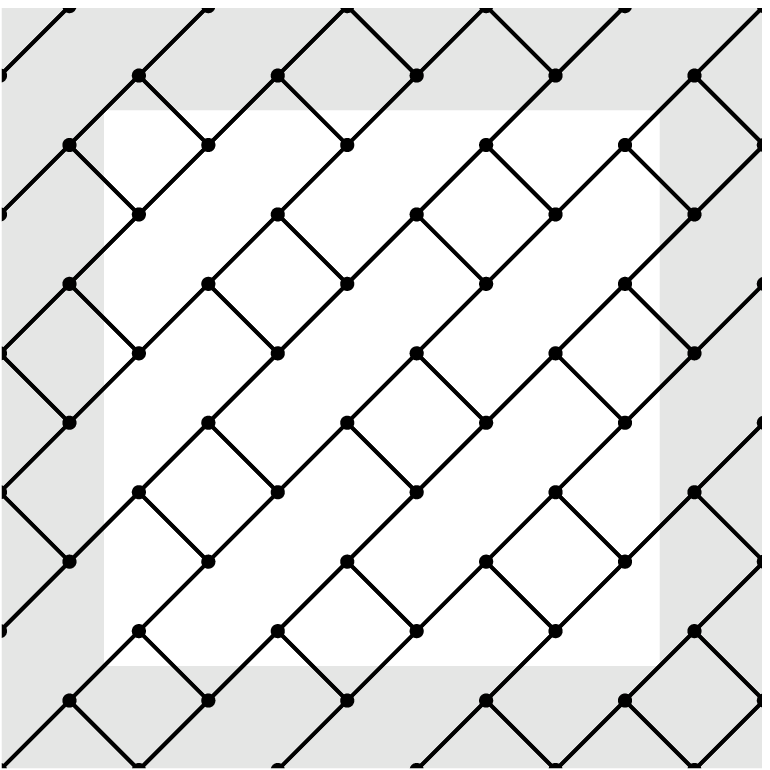}} }}}
\hskip 1cm
\rotatebox{-90}{\reflectbox{\rotatebox{90}{\resizebox{!}{100pt}{\includegraphics{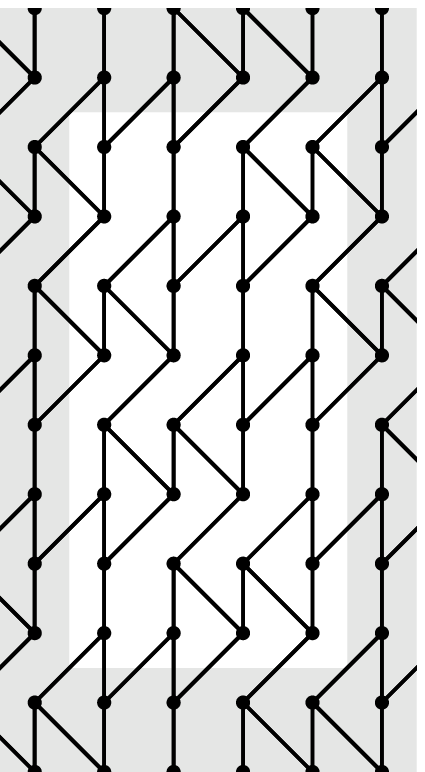}}}}}
\]

\caption{The lattice $(4,8^2)$ drawn with vertices a subset of $\Z^2$,
in original form and in squashed form.}
\label{fig_squash}
\end{figure}

The first step is to transform each lattice so that translations through
some small constant $C$ in the $x$- and $y$-directions act as isomorphisms.
Such a representation of the lattice $(4,8^2)$ with $C=8$ 
is shown
on the left of Figure~\ref{fig_squash}; in this drawing, the vertex
set consists of all points $(x,y)\in \Z^2$ with $x+y$ odd (for some of
the lattices we use $x+y$ even).
The white central portion of the figure shows a square region $S_u$
with scale parameter $s=8$. It is this drawing of the lattice that we consider
when defining $S_u$, $R_e$ and $E_e$.

Note that $s$ must be a multiple of $C$, so that all squares $S_u$ induce isomorphic
subgraphs of the lattice. For the lattices with mirror symmetry (all except
for $(3^4,6)$ and its bond- and site- duals), we choose a representation
with the line $x=y$ or $x=-y$ as an axis of symmetry; a rectangle $R_e$
corresponding to a horizontal bond $e$ may be mapped into a rectangle
$R_f$ corresponding to a vertical bond by a reflection in either of these lines,
so this ensures that all rectangles $R_e$ induce isomorphic subgraphs of $\La$;
thus our program need only evaluate $\Pr_p(E_e)$ for one fixed (horizontal) bond of $\Z^2$.
For the lattices without such symmetry, we run the same program on two
drawings of the lattice, related by reflection in the line $x=y$; the horizontal
rectangle considered for the second drawing corresponds to a vertical one
in the first.

Most of the representations we use are modifications of those shown in
Figure 3 of Suding and Ziff~\cite{SudingZiff}, most of which have a
horizontal axis of mirror symmetry. Since we want a diagonal axis here,
we have rotated may of the representations by 45 degrees,
obtaining a graph on points in $\Z^2$ with $x+y$ even. Our
representations for all 11 Archimedean lattices $\La$ and their planar
duals $\La^\star$ are shown in the appendix. We have omitted any
lattice for which the critical probability is known exactly.
(In each case the 
site dual $\Lax$ is represented in the same way as $\La$, but
with additional bonds added to every face.)

For those lattices represented with vertices in $\Z^2$ with $x+y$ even,
it is computationally more efficient to modify the representation to make
it more compact, by mapping  $(x,y)$ to $(\lfloor x/2\rfloor,y)$, say. 
An example for the lattice
$(4,8^2)$ is shown on the right of Figure~\ref{fig_squash}. Note that an $s/r$ by
$s$ rectangle in the compact form corresponds   
to an $s$ by $s$ square in the original, where $1/r=1/2$ is the
{\em ratio} by which we have squashed the lattice when compactifying
it. In the program files, this is stored as the field RATIO for each  
lattice; this squashing is undone in the print\_lattice routines.

The program perc.c, available from our website \cite{WEB}
reads in the lattice, assigns states to the sites
or bonds randomly, and then finds the largest open clusters in the left- and right- halves 
$\La_u$ and $\La_v$ of $\La_e$, the subgraph of $\La$ induced by the
sites in $R_e$. (In fact, to avoid using too much memory, these two
processes are done concurrently, see below for details).
Finally, it tests whether these open clusters are joined in $\La_e$.
The open clusters are found using a simple incremental algorithm
that scans $S_u$ from the left and $S_v$ from the right.
The method used to find the largest open cluster is (a simplified form of)
that of Balister, Bollob\'as and Walters~\cite{BBWsquare}, and works as follows. We
divide the square $S_u$ into strips which are narrow  but are
sufficiently wide that no edge jumps an entire strip (i.e., all the
edges meeting the strip are entirely contained in the union of the
strip and its two neighbouring strips).

We find the component structure of the open subgraph restricted to the two
left most strips. Then we look at the next strip and find the new
component structure formed. At each stage we have an equivalence
relation on the vertices in a strip where two vertices are equivalent
if they are in the same open cluster in the part of $S_u$ to the left of
the current strip. We also keep track of the size of each of these
clusters, and the size of the largest open cluster we have seen so
far. When we get to the right hand edge of the square $S_u$ we know exactly
which vertices (if any) in that strip are part of the largest
open cluster of $\La_u$.

We repeat the process on $S_v$ but working from right to left. Finally
we add the edges between the right most strip of $S_u$ and the left
most of $S_v$ and see whether the largest open clusters in each are
joined.

The important thing to note about this algorithm is that the storage
required is proportional to the side length of $S_u$, i.e. to $s$, not to
the area of $S_u$.

\section{Results}

For our percolation bounds see Table~\ref{table1}. For full results,
including numbers of successes, please see our website \cite{WEB}.
Note that in the 400 simulations associated with each bound (or with
each attempt to obtain a bound) we have always seeded the random
number generator with 400 consecutive seeds starting from
$12345678$. This means that the exact results of our simulations
should be reproducible as a way of checking the program. Also, it
shows that we have not performed many different runs and finally
chosen seeds that work!

The computations were performed running in the background on
around 70 (mostly fairly old) computers in the Department of
Pure Mathematics and Mathematical Statistics, University of Cambridge, over a period
of around 2 weeks. This was made much easier by the fact that the department uses
Linux rather than Windows! We are grateful to the computer officer, Andrew Aitchison,
for technical assistance.

\begin{table*}
\begin{center}
\begin{tabular}{|l||c|c||c|c||}\hline
Lattice & Site  & Width  & Bond  & Width\\\hline
Square     & [0.5925,0.5930]	&$5\times10^{-4}$	&0.5	&	0	\\
Triangular & 0.5	                &0	                & $2\sin(\pi/18)$	&	0	\\
Hexagonal  &[0.6968,0.6973]	&$5\times10^{-4}$	& $1-2\sin(\pi/18)$	&	0	\\
Kagom\'e     &$1-2\sin(\pi/18)$		&0	&[0.52415,0.52465]	&$5\times10^{-4}$ \\
$(3,12^2)$    &$\sqrt{1-2\sin(\pi/18)}$		&0	&[0.7402,0.7407]	&$5\times10^{-4}$\\
$(3,4,6,4)$    &[0.6216,0.6221]	&$5\times10^{-4}$	&[0.5246,,0.5251]	&$5\times10^{-4}$\\
$(3^3,4^2)$      &[0.5500,0.5505]	&$5\times10^{-4}$	&[0.4194,0.4199] 	&$5\times10^{-4}$\\
$(3^2,4,3,4)$      &[0.5506,0.55105]	&$4.5\times10^{-4}$	&[0.4139,0.4144]	&$5\times10^{-4}$\\
$(3^4,6)$      &[0.57925,0.57975]   &$5\times10^{-4}$	&[0.4341,0.4345]	&$4\times10^{-4}$\\
$(4,6,12)$     &[0.7476,0.7480]	&$4\times10^{-4}$	&[0.6935,0.6940]	&$5\times10^{-4}$\\
$(4,8^2)$        &[0.7295,0.7300]	&$5\times10^{-4}$	&[0.6766,0.6770]	&$4\times10^{-4}$\\\hline
\end{tabular}
\end{center}
\caption{Rigorous 99.9999\% confidence intervals for critical probabilities for site and bond percolation}\label{table1}
\end{table*}

\section{Conclusions}\label{concl}

We have shown that it is practical to use the method of Balister, Bollob\'as and Walters~\cite{BBWsquare}
to obtain narrow confidence intervals for the critical probabilities for site
and bond percolation on Archimedean
lattices. Unlike the (presumably) much more precise estimates obtained by other methods,
these intervals come with mathematically guaranteed error bounds. The intervals
are much narrower than those that can be 100\% proved, and the error probabilities
are very small; the running time does not increase much with a large decrease in the
desired error probability, so a probability that is in practice zero (here 1 in a million for each lattice)
may be achieved.

We have tried to keep the computations relatively simple; there is no point in using an algorithm
that is proved correct if it is not possible to verify the computer program used.
At the cost of more complicated programing, better results could be obtained in two ways.
Firstly, the current program could be made to cache better and hence run faster
by scanning the rectangle $R_e$ in a more complicated manner: this $2s$ by $s$ rectangle
could be broken down into $k$ by $k$ squares small enough that the boundary of one
square fits into the processors primary cache, and these squares could then be processed
column by column. The overall storage requirement is approximately the same (one entire
column must be stored), but the frequency of cache misses is reduced by a factor of about~$k$.

A more significant improvement could be obtained by considering a different event $E_e$:
let $E_e$ be the event that there is an open path crossing $R_e$ from left to right, and
that there is an open path crossing the left-hand end square of $R_e$ from top to bottom.
As noted by Bollob\'as and Riordan~\cite{BRKesten}, 
for example, this event still has the property \eqref{Acond}. (Essentially this observation
was used by Balister, Bollob\'as and Walters~\cite{BBWsquare} in obtaining a lower
bound on the critical parameter for a certain continuum percolation model.)
Also, the scaling behaviour of $\Pr_p(E_e)$ near $\pc$ should be the same
as for the event considered here. The gain is that whether or not $E_e$ holds in a given
configuration can be tested faster, using an interface
following algorithm of the type used by Ziff and Cummings
\cite{ZiffCummings} in 1984,
for example. Assuming conformal invariance,
the expected length of the interface is $s^{2-\alpha_3}=s^{4/3}$, where $\alpha_3=(3^2-1)/12$
is the multi-chromatic $3$-arm exponent for $SLE_6$.
Note, however, that to use this algorithm in practice without running into memory/caching problems,
one needs to generate the state of each site/bond from a pseudo-random {\em function},
rather than a pseudo-random number generator.

All the exactly known critical probabilities associated to Archimedean
lattices are roots of (simple) polynomial equations with integer
coefficients.  While it is easy to construct other lattices whose
critical probabilities may be found in this form (see Ziff
\cite{ZiffCellDualCell}, for example), it may well be that there are
no such expressions for the remaining Archimedean lattices, although
some have been conjectured. In particular, Wu \cite{Wu} conjectured
that for bond percolation on the Kagom\'e lattice, $\pcb=0.524429$,
a root of the equation $p^6-6p^5+12p^4-6p^3-3p^2+1=0$.  Tsallis
\cite{Tsallis} conjectured the values $\pcb=0.522372 $ and
$\pcb=0.739830$ for bond percolation on the Kagom\'e and $(3,12^2)$ lattices,
respectively. Tsallis's conjectures have been effectively ruled out
some time ago by experimental estimates (see \cite{Yonezawa}, for
example); they are rather far from the current best estimates of
$0.5244053$ and $0.74042195$. They have not yet
been rigorously disproved, although the latest results of May and
Wierman \cite{MayWierman} come close.  For both lattices, our results
provide a rigorous `99.9999\% disproof' of Tsallis's values - they lie
outside our rigorous 99.9999\% confidence intervals.

\textheight=625pt

Wu's conjectured value for $\pcb(K)$ seems to be much closer to the
truth; it is well within the confidence interval we
obtain. Nevertheless, it is still believed to be false; see Ziff and
Suding \cite{ZiffSuding97}, for example. More recently, Scullard and Ziff
\cite{ScullardZiff} have predicted certain values for $\pcb$ for the Kagom\'e and
$(3,12^2)$ lattices, using a heuristic version of the star-triangle
transformation. Although they leave open the `possibility' that one of
these values might be exact, there seems no reason (to us, or,
apparently, to them) to really believe this: the method is (as they
admit) non-rigorous, and the value obtained in the same way for the
Kagom\'e lattice (given earlier by Hori and Kitahara without
derivation) is outside the error bounds of existing experimental
results.

\section{Appendix}

\begin{figure*}
\[
\begin{array}{c}
\begin{array}{cccc}
\epsfig{file=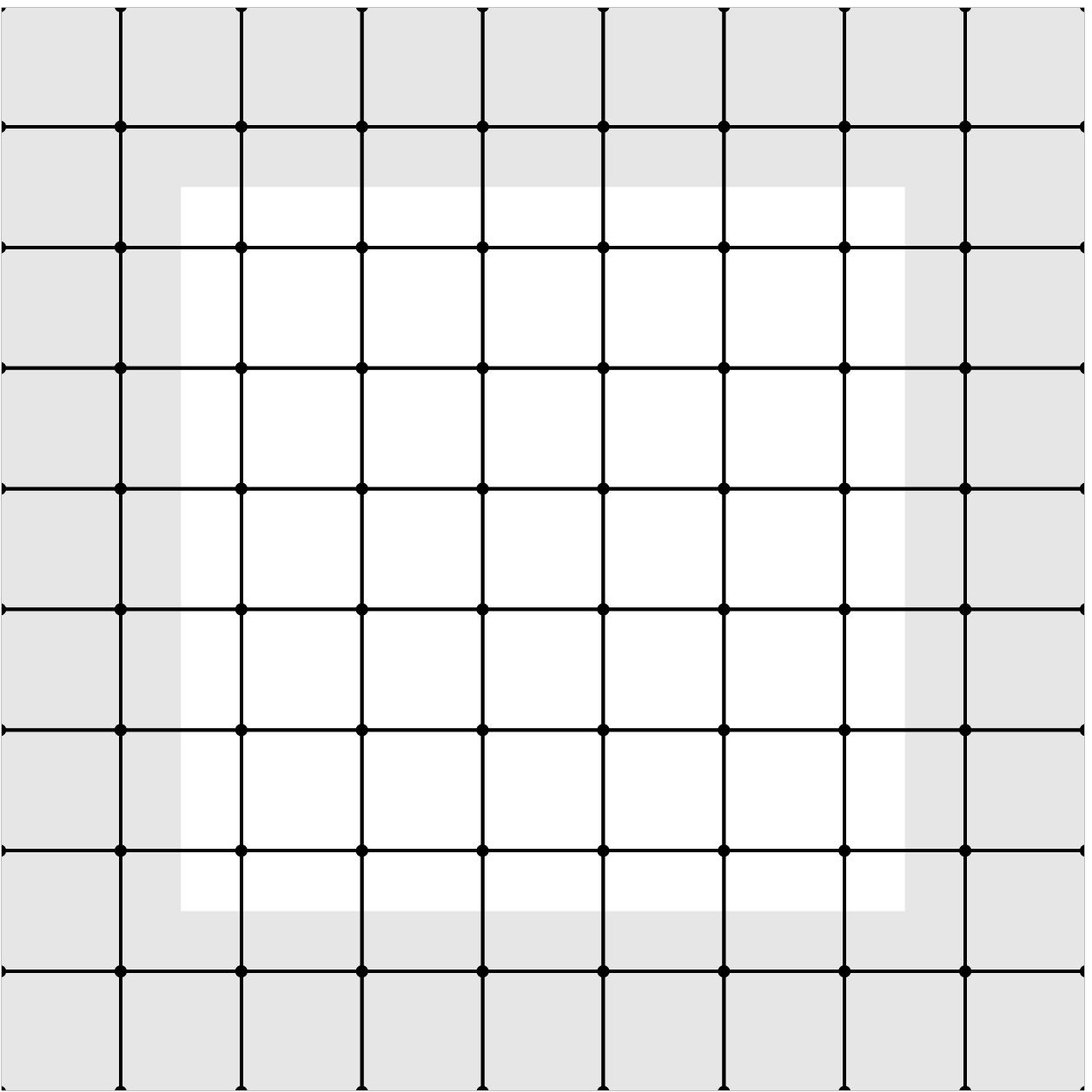,height=80pt} &
\epsfig{file=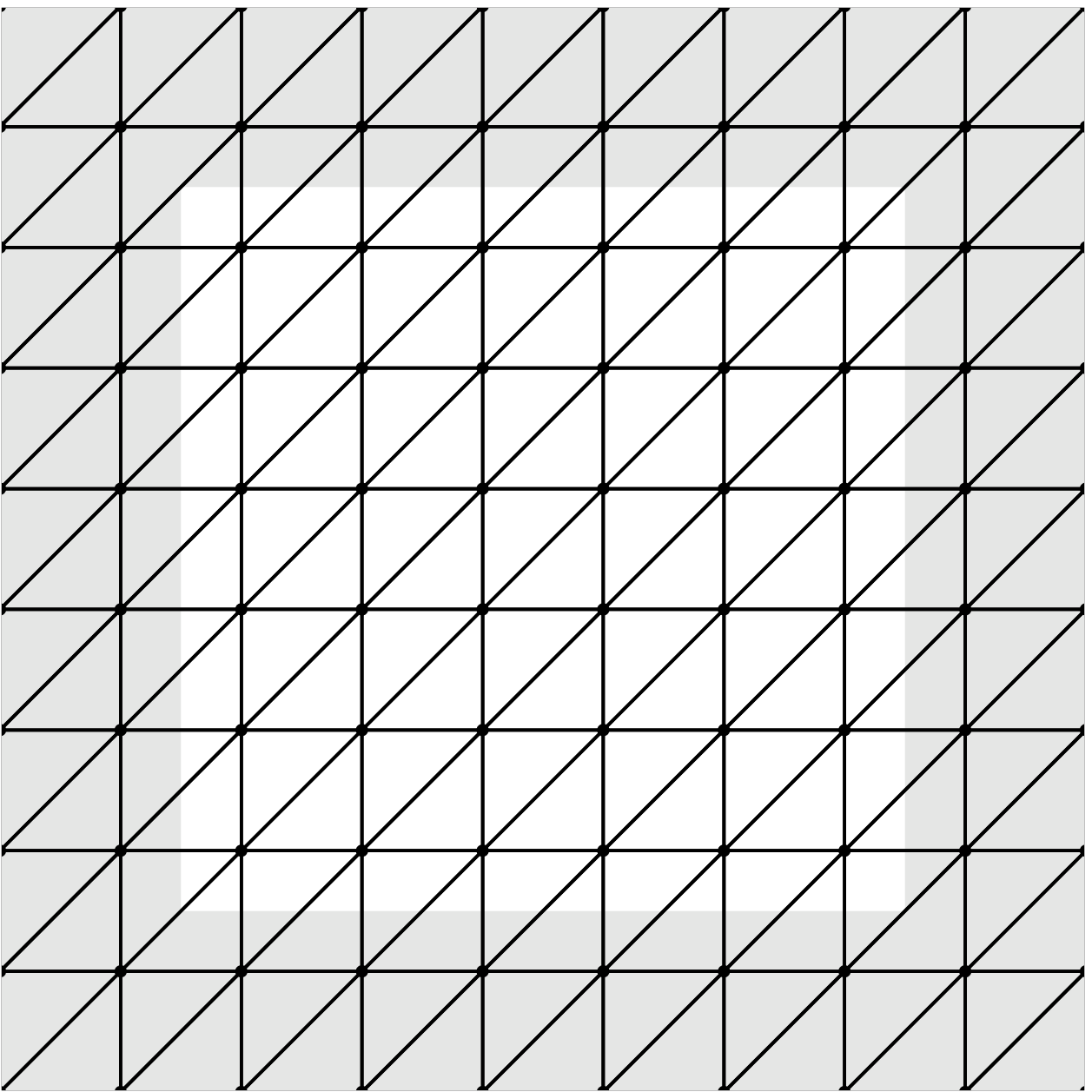,height=80pt} & 
\epsfig{file=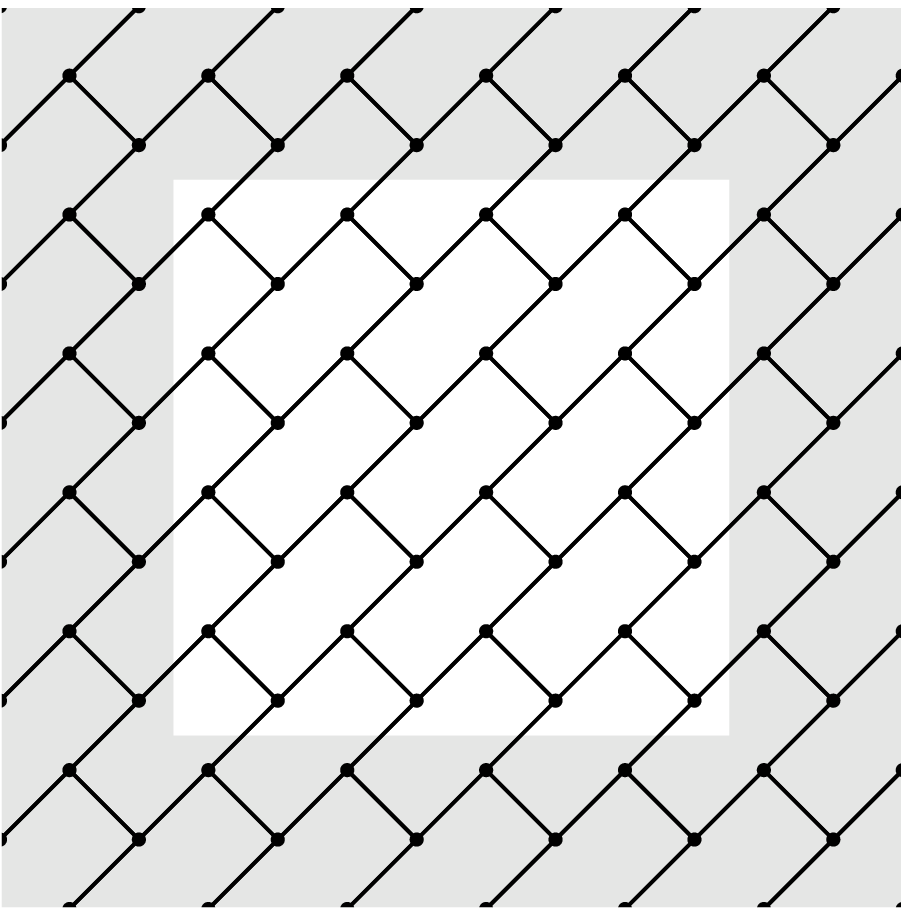,height=80pt} &
\epsfig{file=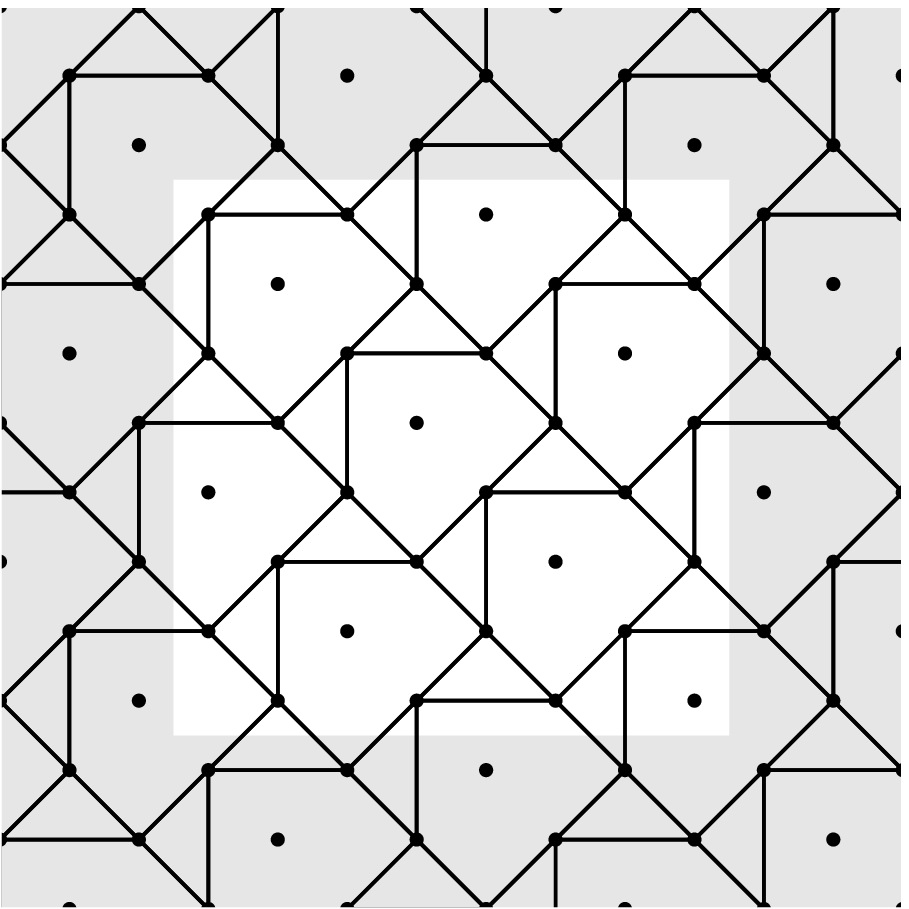,height=80pt} \\
 \hbox{Square: }(4^4) & \hbox{Triangular: }(3^6) & \hbox{Hexagonal:
  }(6^3) & \hbox{Kagom\'e: }(3,6,3,6) \\ \\
\end{array}\\
\begin{array}{ccc}
\epsfig{file=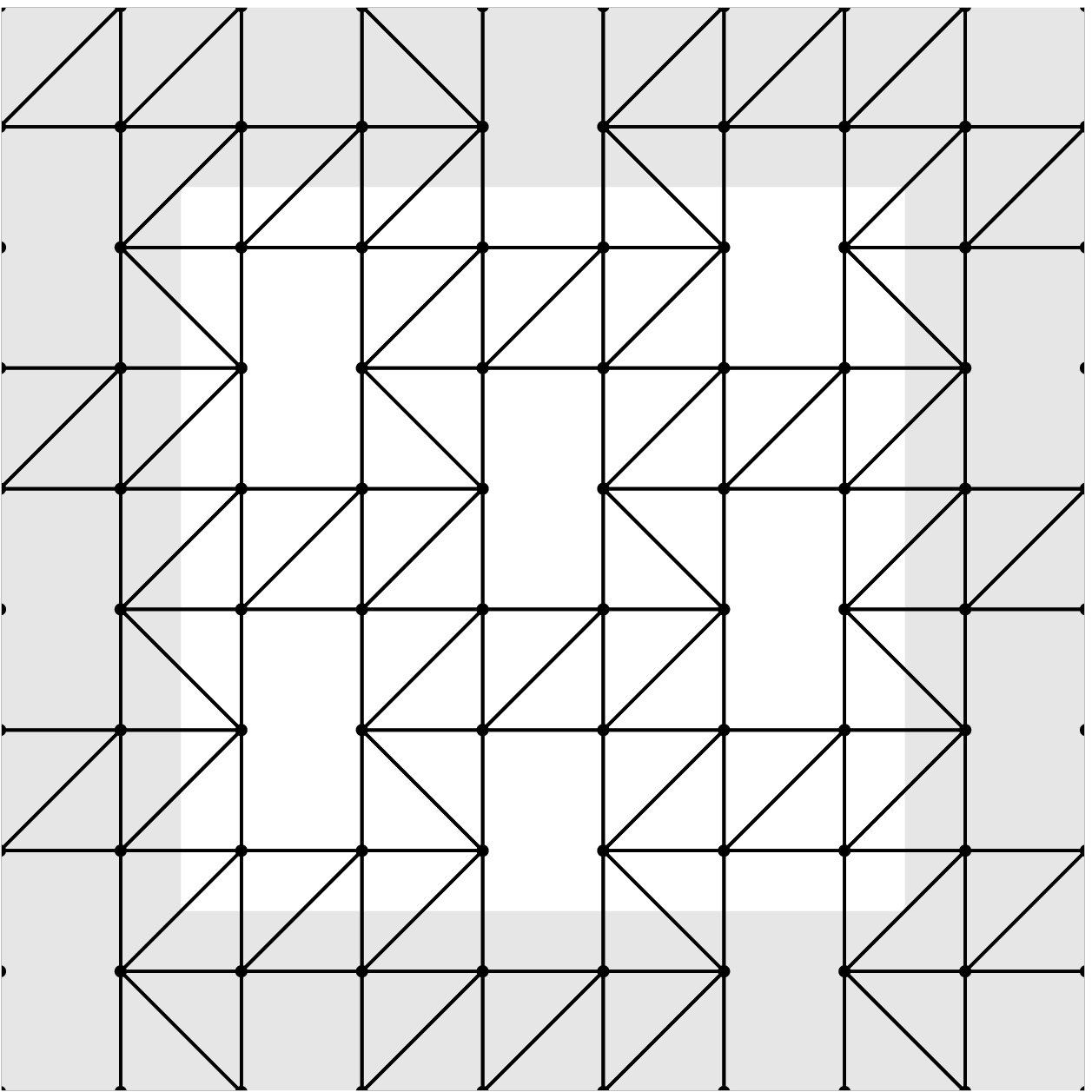,height=80pt} &
\epsfig{file=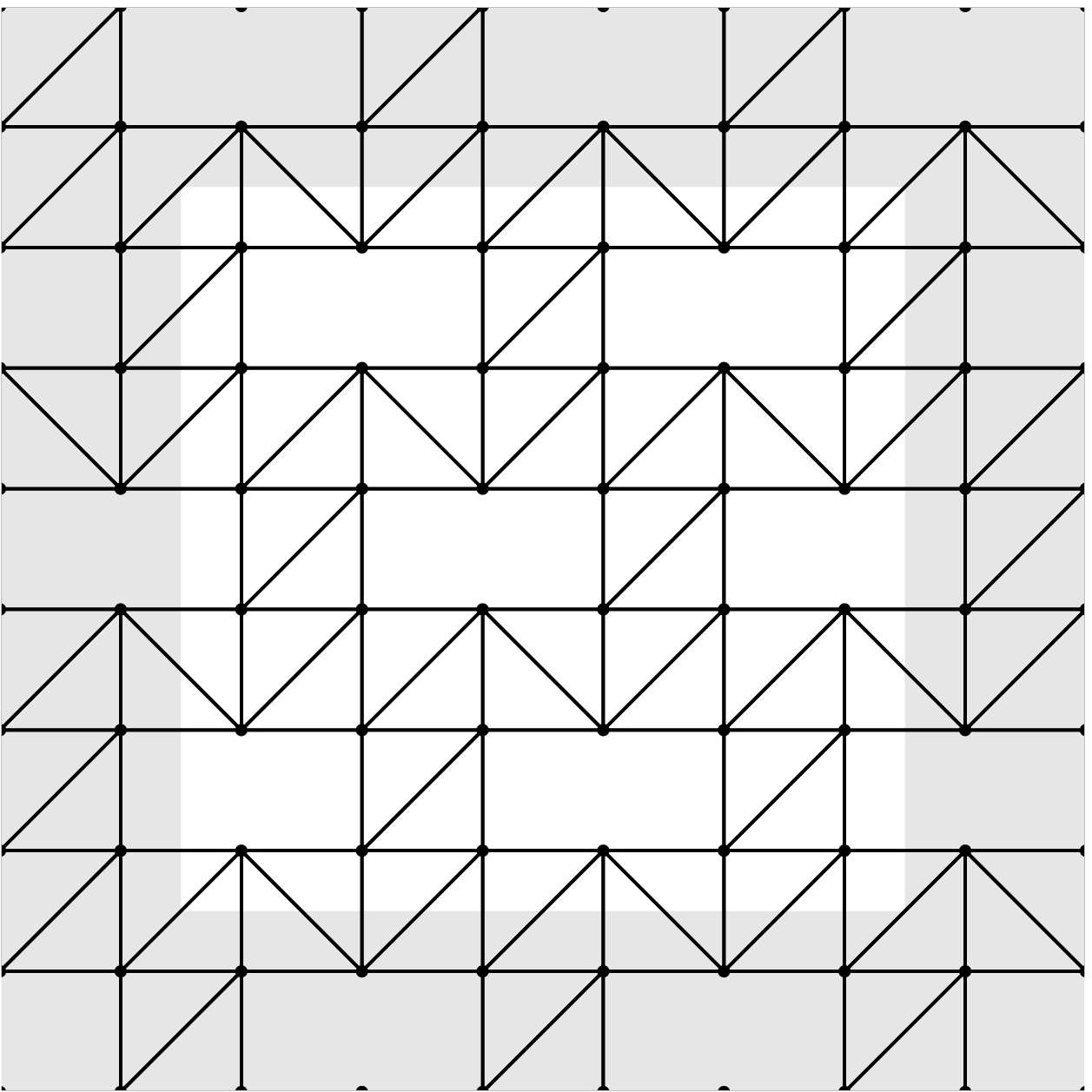,height=80pt} &
\epsfig{file=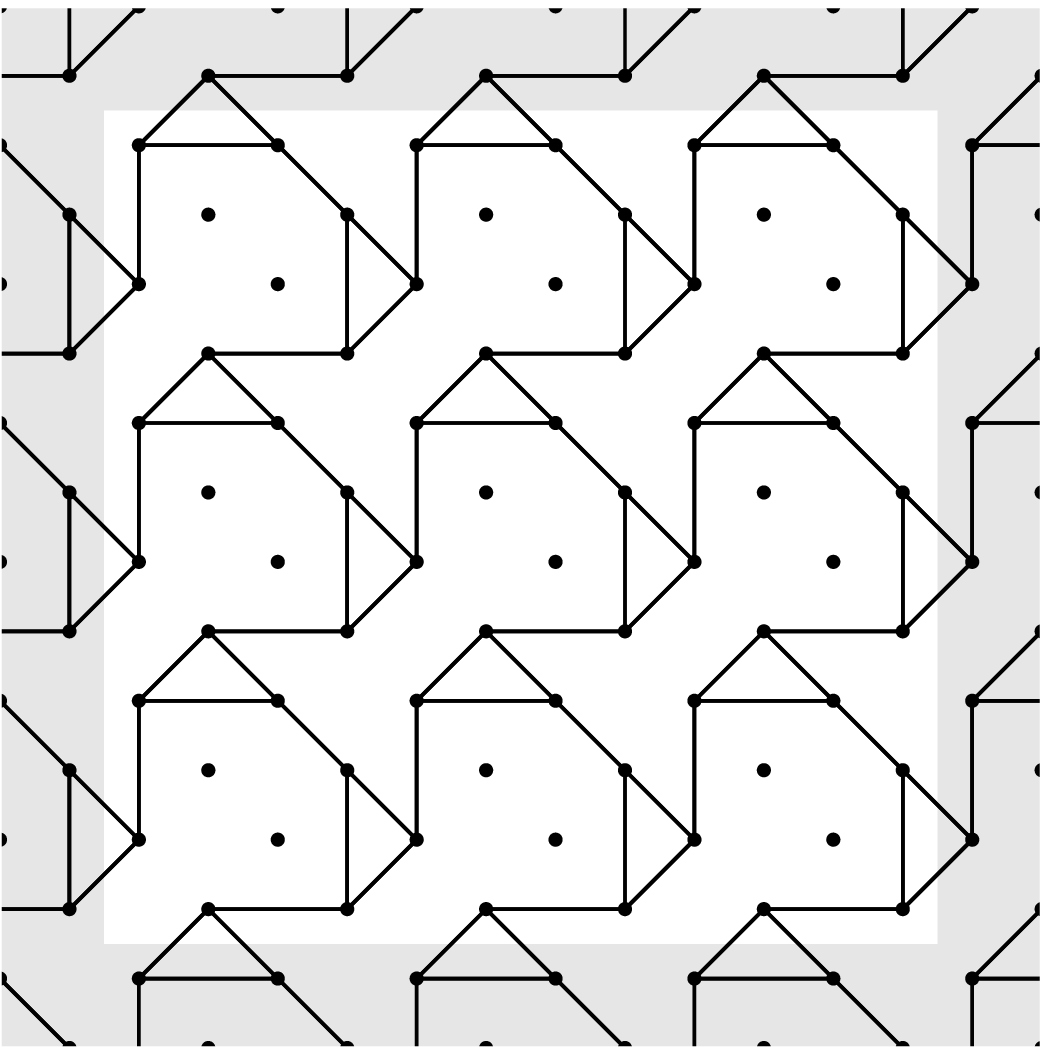,height=80pt}\\
\hbox{$(3^4,6)$ Horizontal} & \hbox{$(3^4,6)$ Vertical} & (3,12^2)   \\ \\
\end{array}\\
\begin{array}{ccc}
\epsfig{file=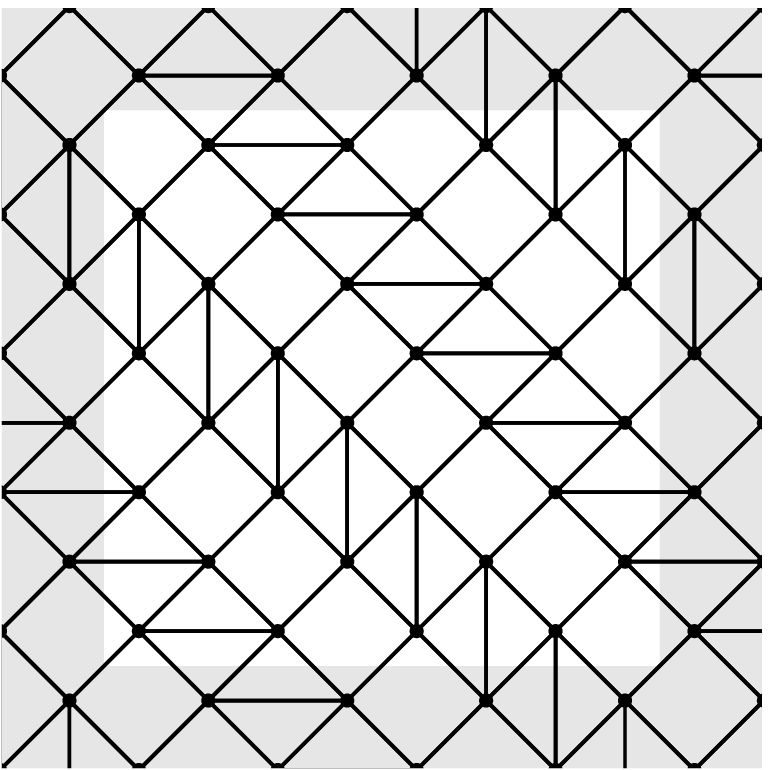,height=80pt} &
\epsfig{file=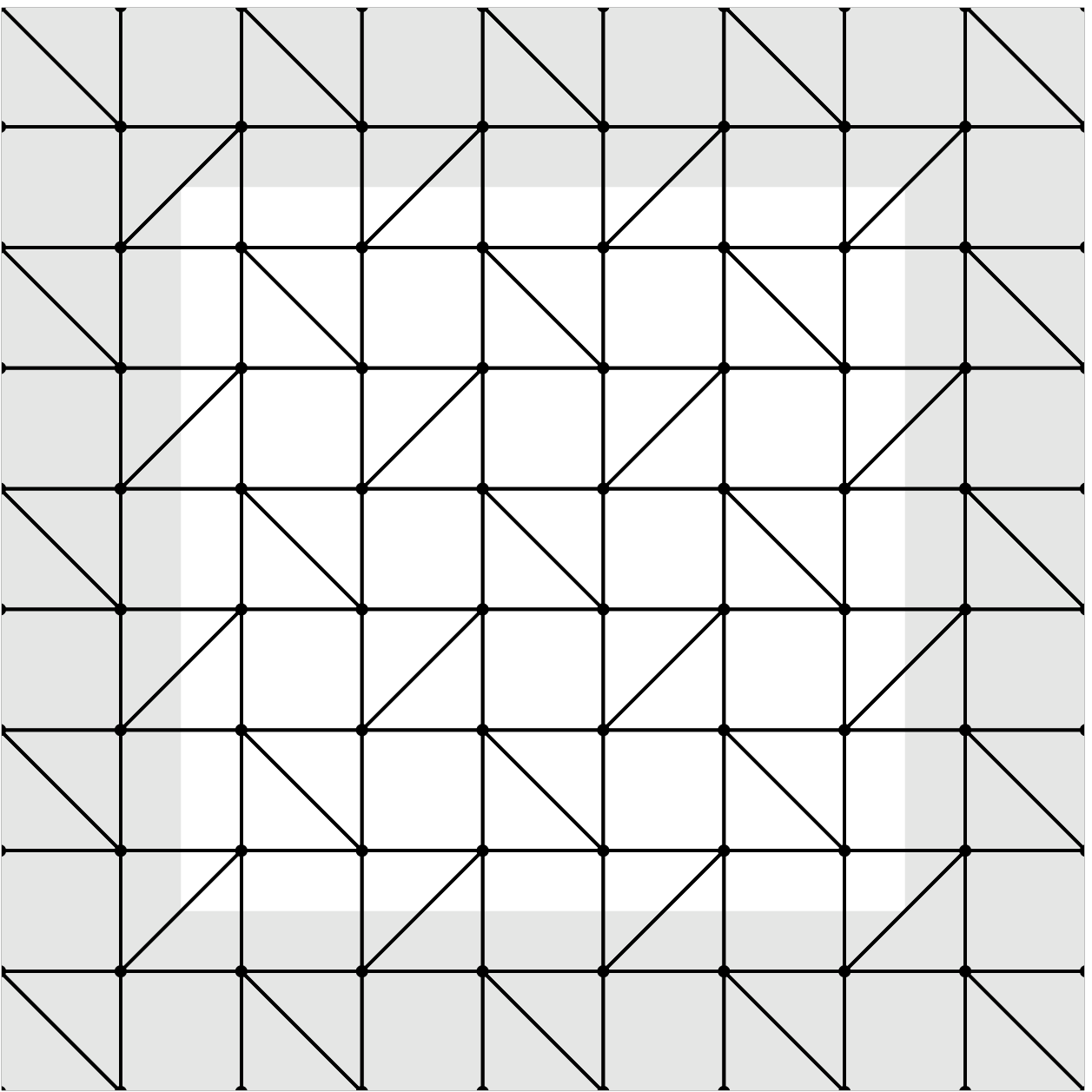,height=80pt} &
\epsfig{file=488.eps,height=80pt} \\
 (3^3,4^2) & (3^2,4,3,4) & (4,8^2) \\ \\
\end{array}\\
\begin{array}{cc}
\epsfig{file=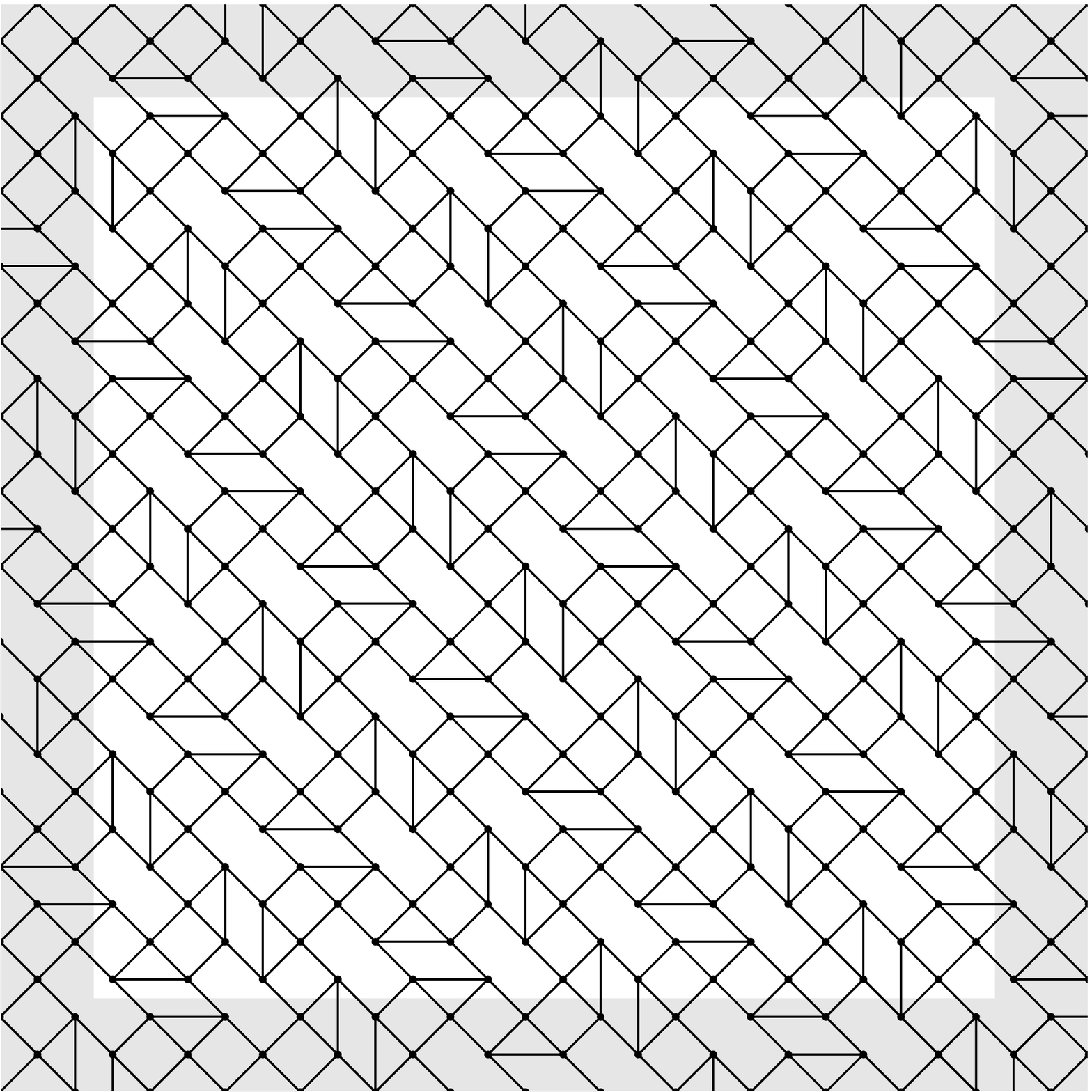,height=150pt} &
\epsfig{file=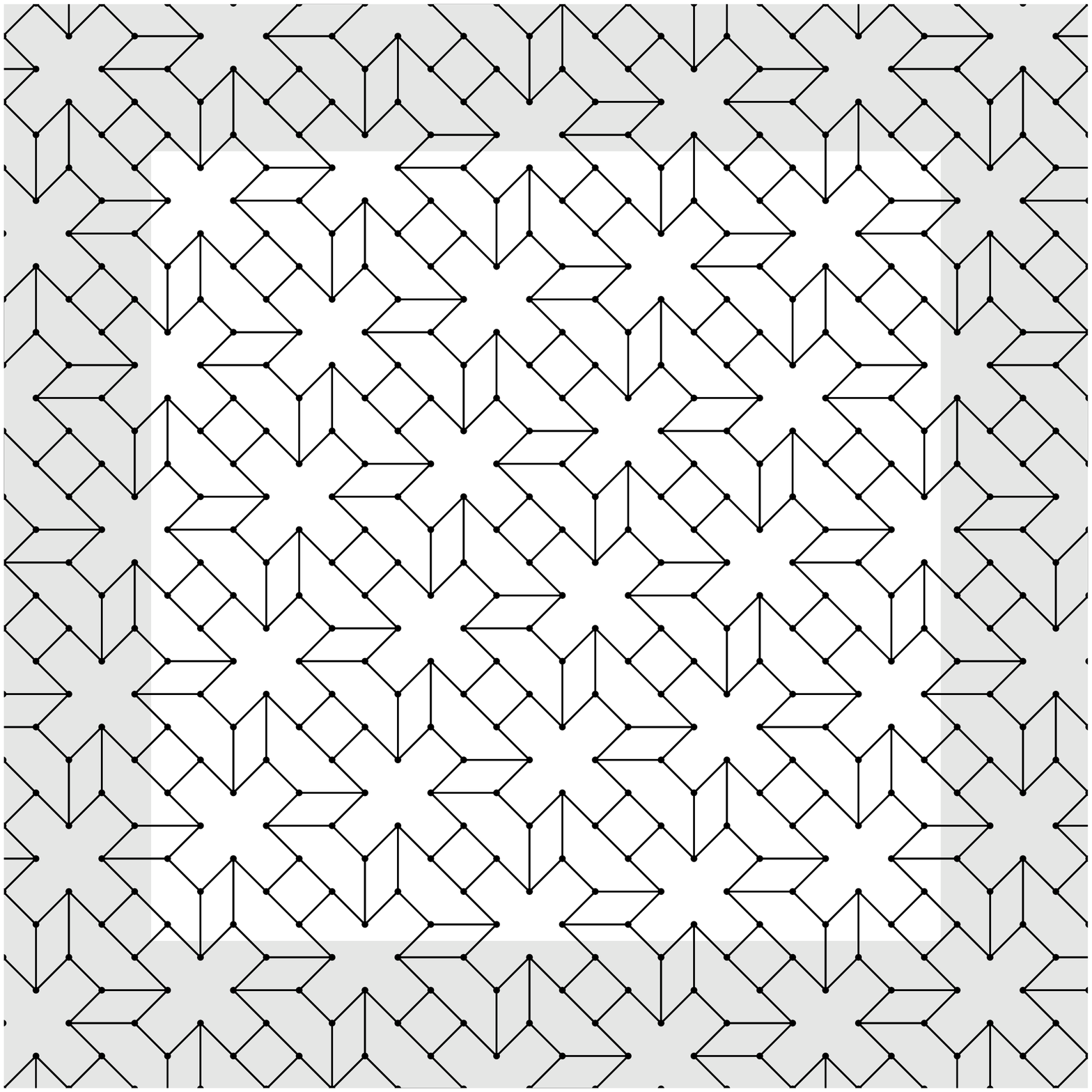,height=150pt} \\
(3,4,6,4)& (4,6,12)\\
\end{array}\\
\end{array}
\]
\caption{The 11 Archimedean lattices  on a square grid.}\label{fig_sqlat}
\end{figure*}

\begin{figure*}
\[
\begin{array}{c}
\begin{array}{ccc}
\resizebox{!}{100pt}{\includegraphics{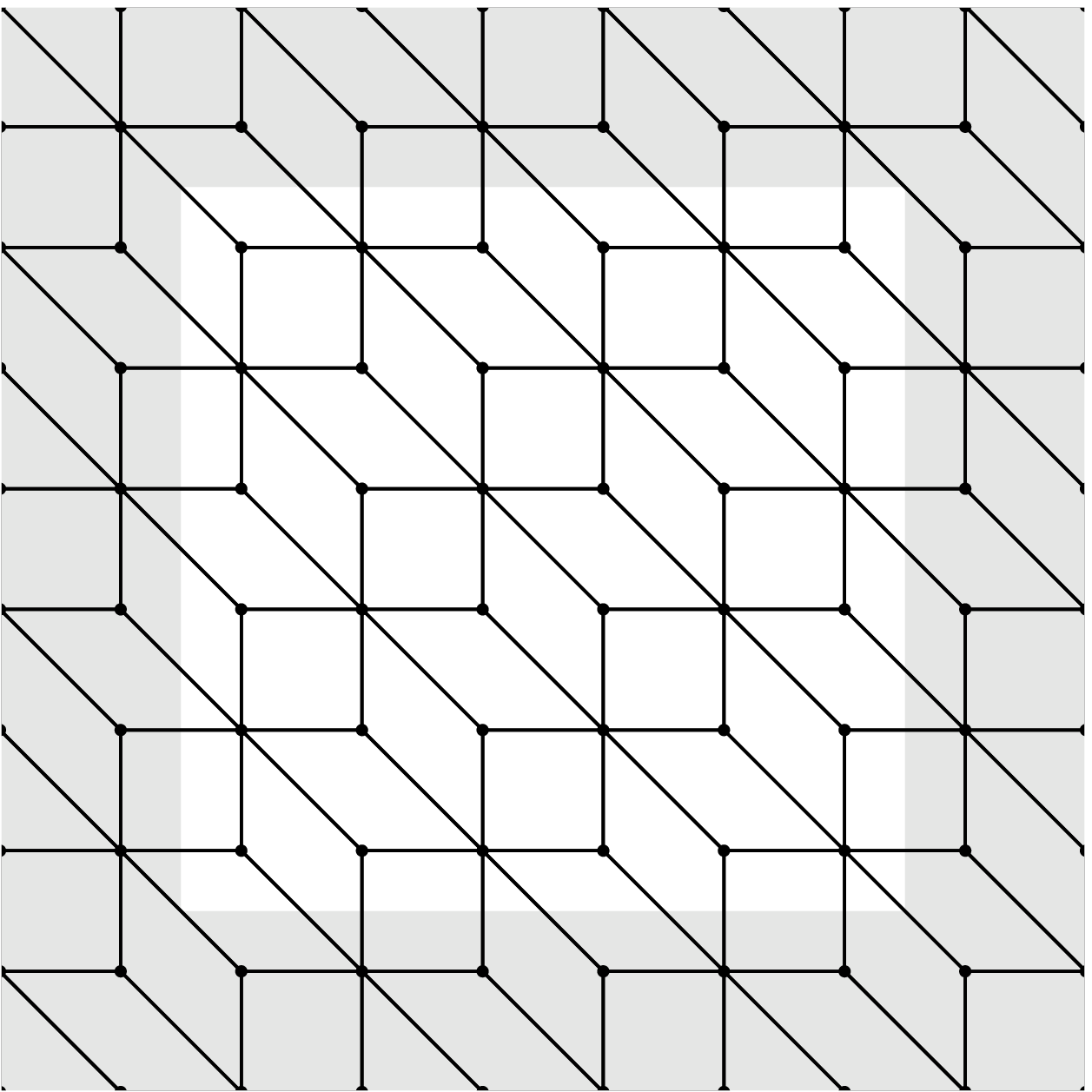}} &
\resizebox{!}{100pt}{\includegraphics{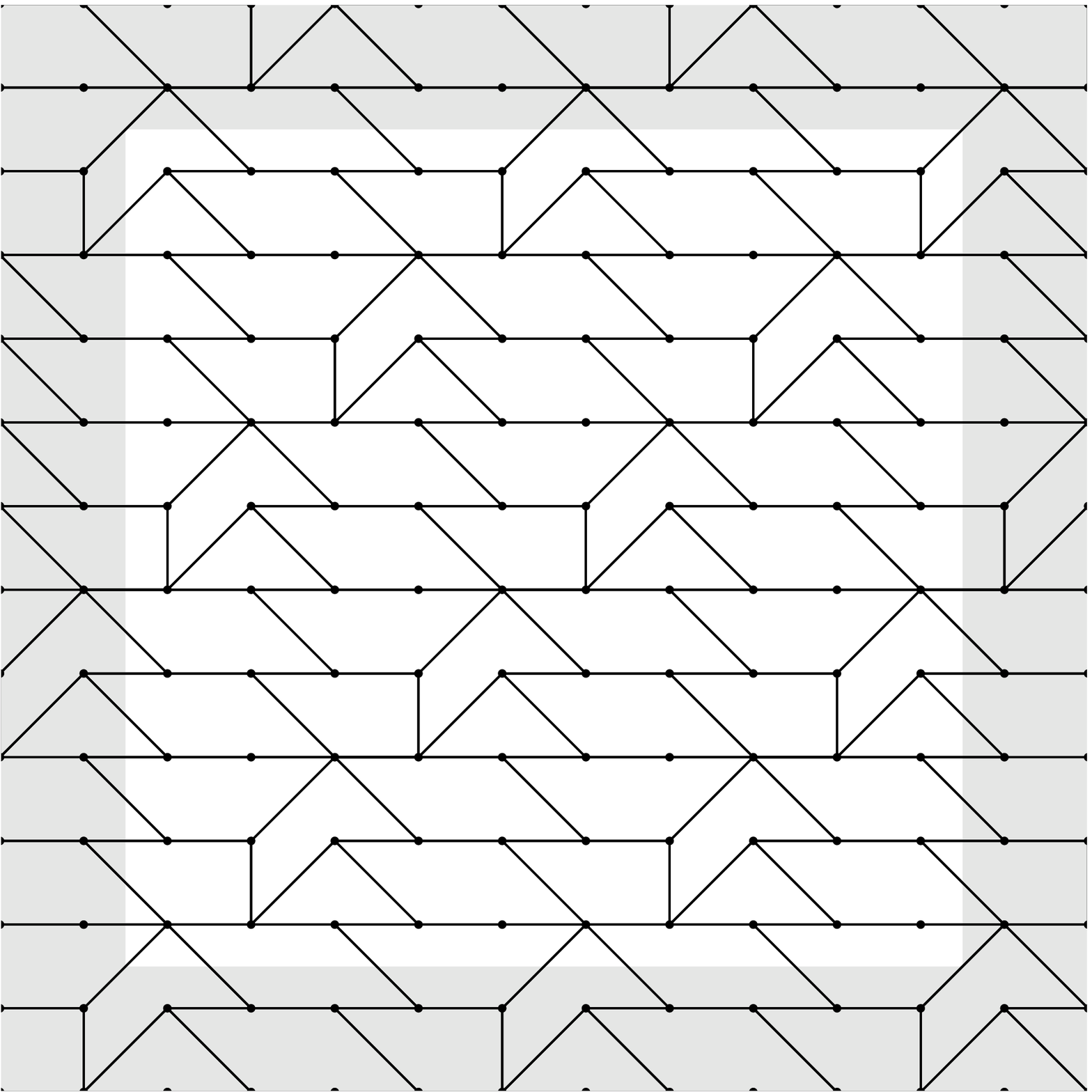}}& 
\resizebox{!}{100pt}{\includegraphics{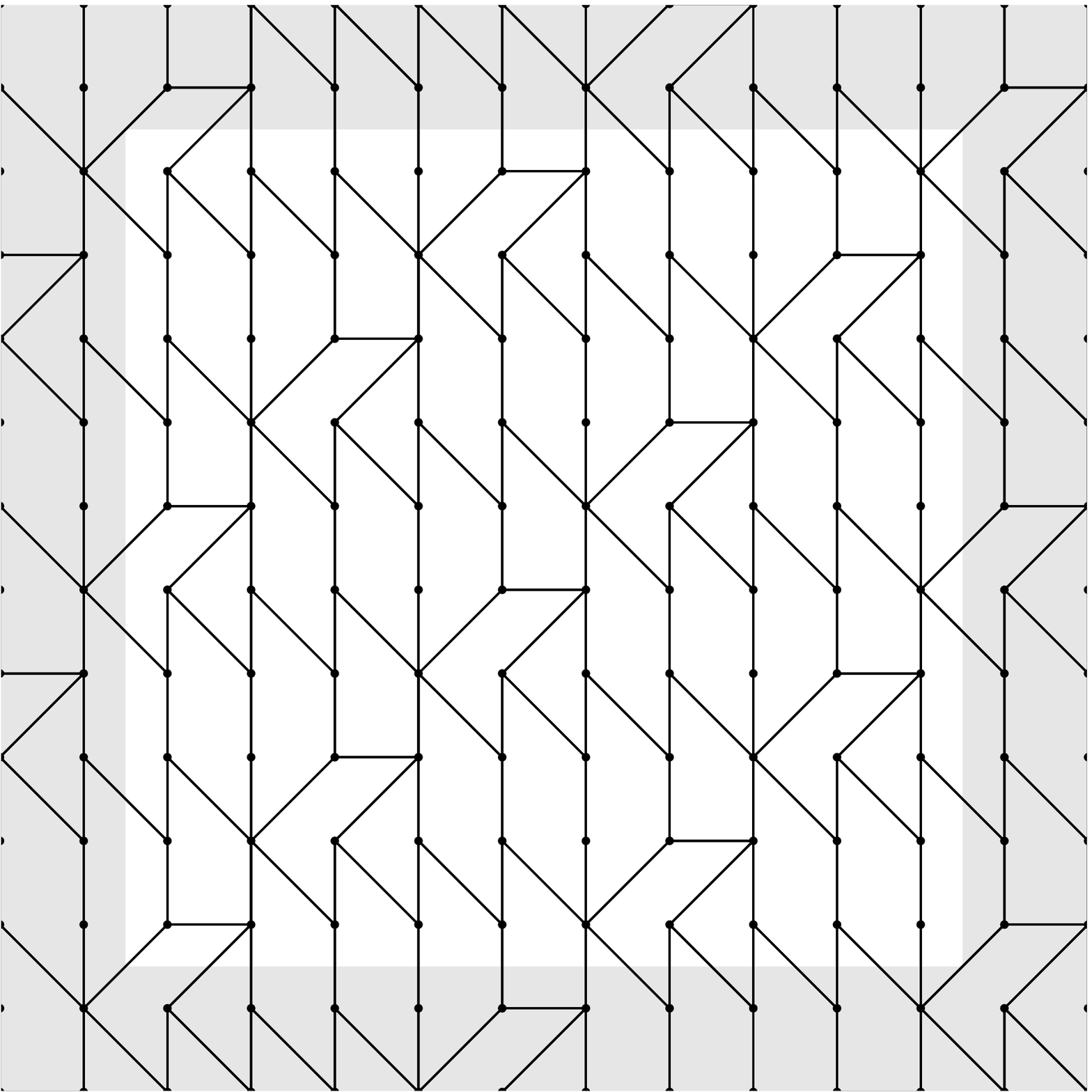}}\\ 
 \hbox{Kagom\'e $(3,6,3,6)$} & \hbox {$(3^4,6)$ Horizontal} & 
\hbox {$(3^4,6)$ Vertical}\\ \\
\end{array}\\
\begin{array}{ccc}
\resizebox{!}{100pt}{\includegraphics{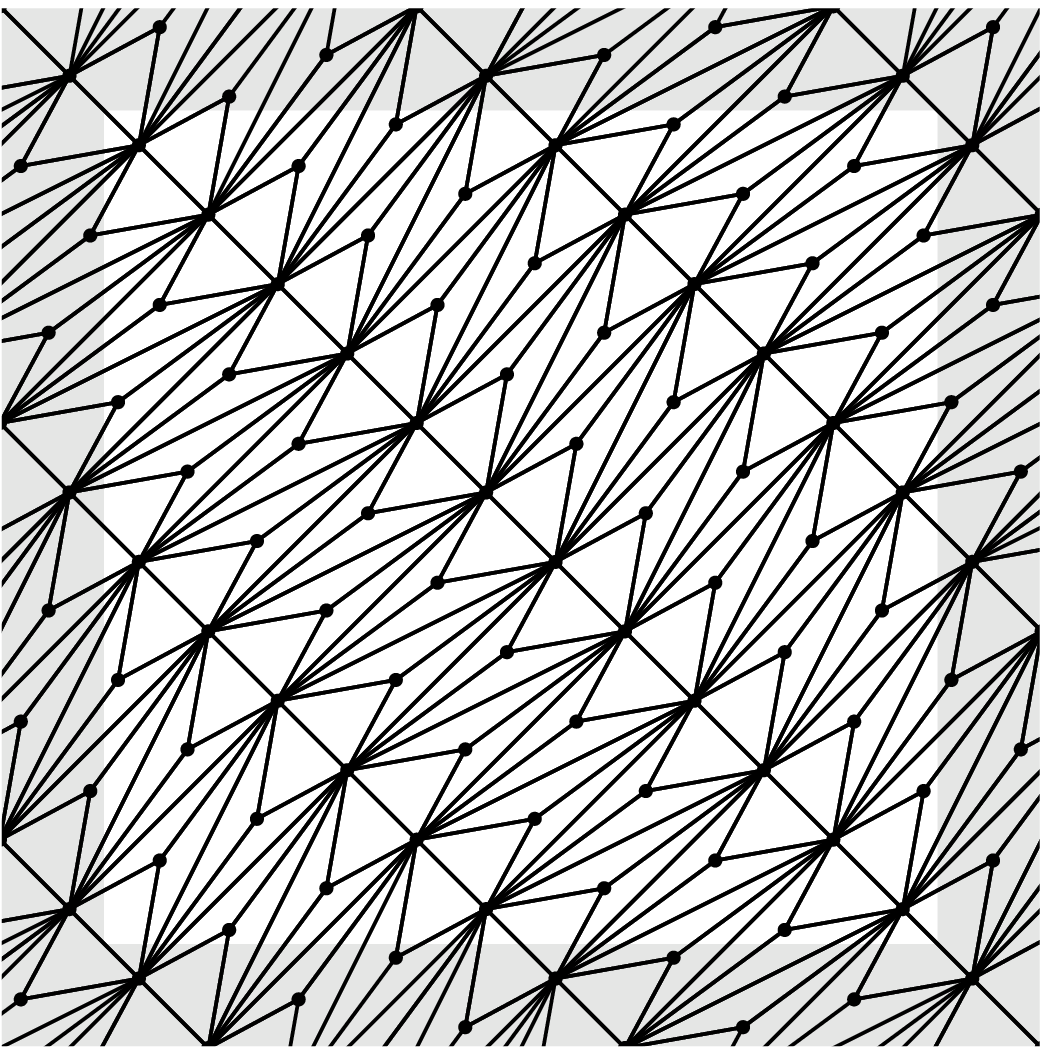}} &
\resizebox{!}{100pt}{\includegraphics{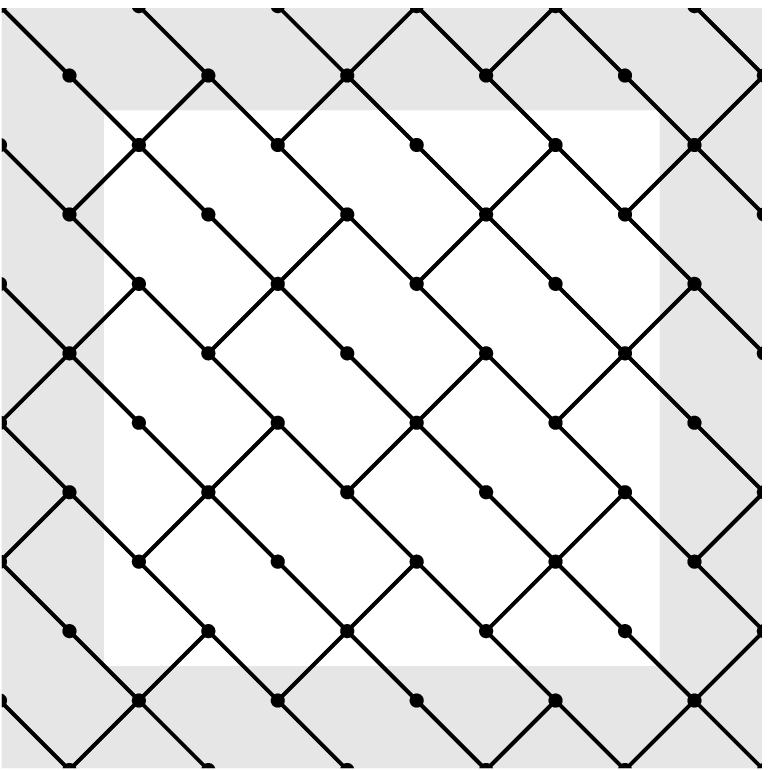}} &
\resizebox{!}{100pt}{\includegraphics{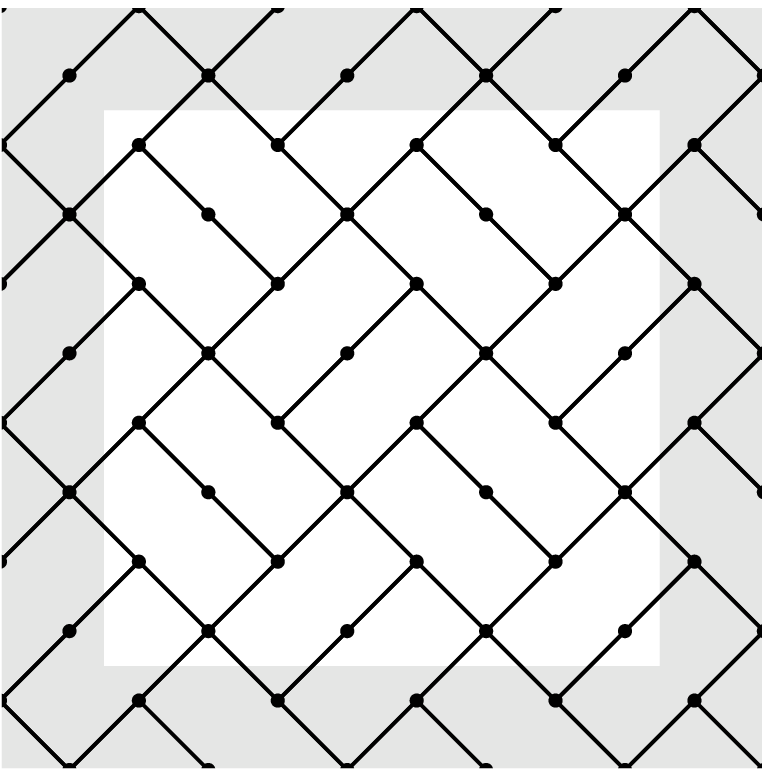}} \\
(3,12^2) & (3^3,4^2) & (3^2,4,3,4)   \\ \\
\end{array}\\
\begin{array}{ccc}
\resizebox{!}{100pt}{\includegraphics{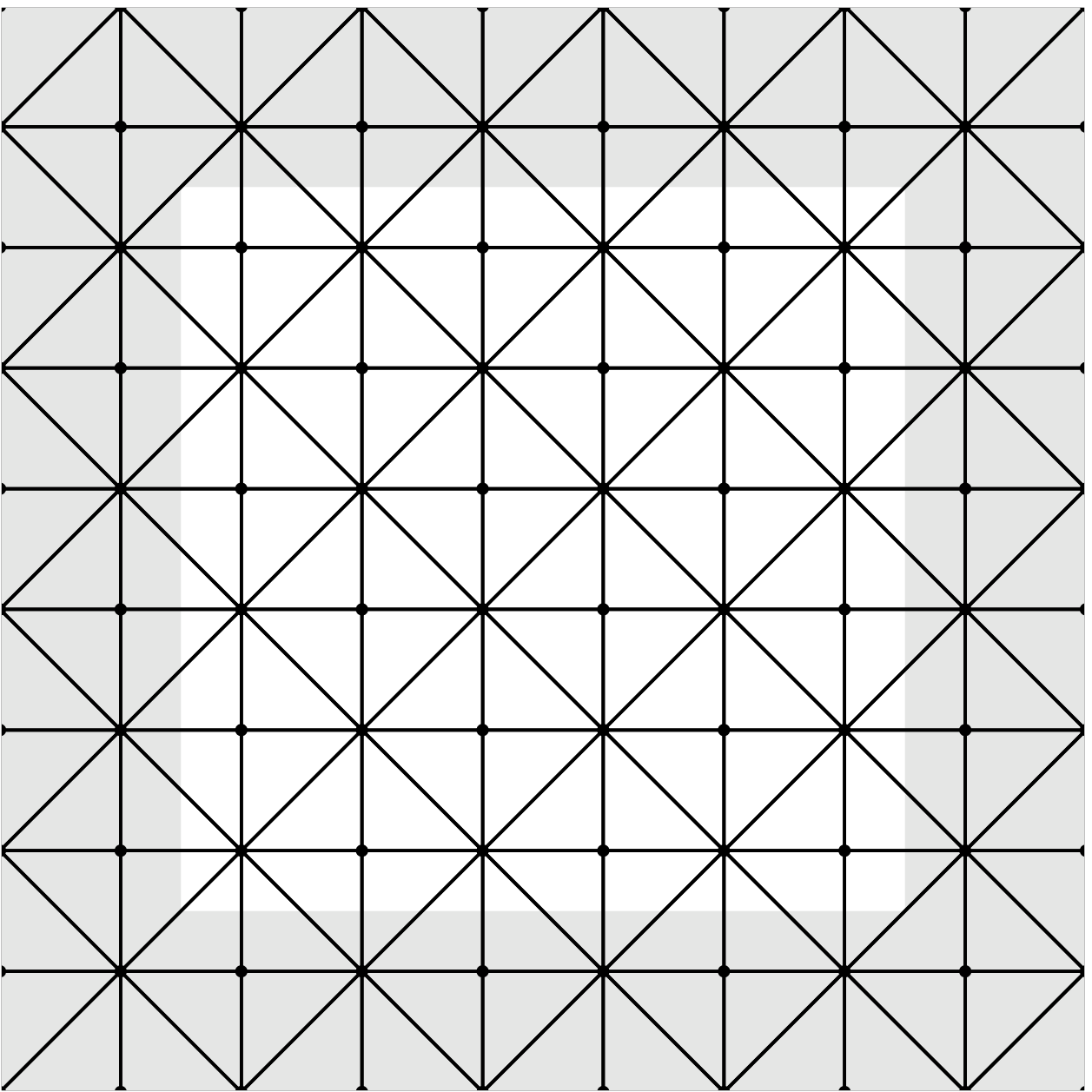}} &
\resizebox{!}{100pt}{\includegraphics{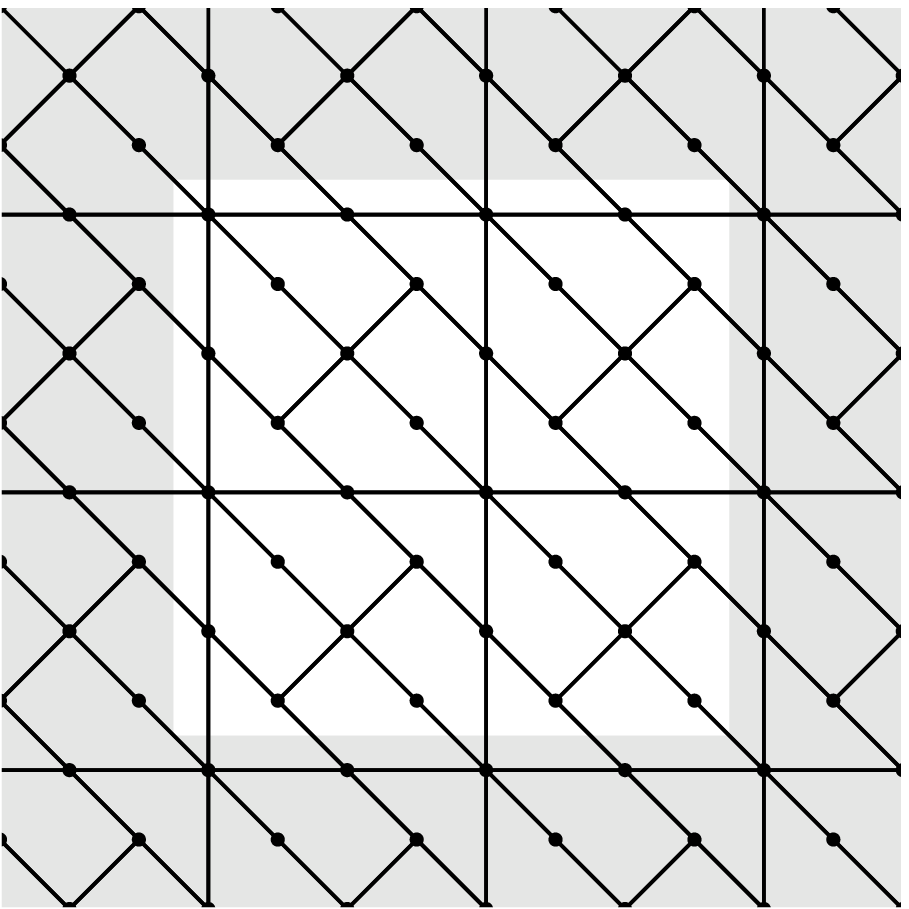}} &
\resizebox{!}{100pt}{\includegraphics{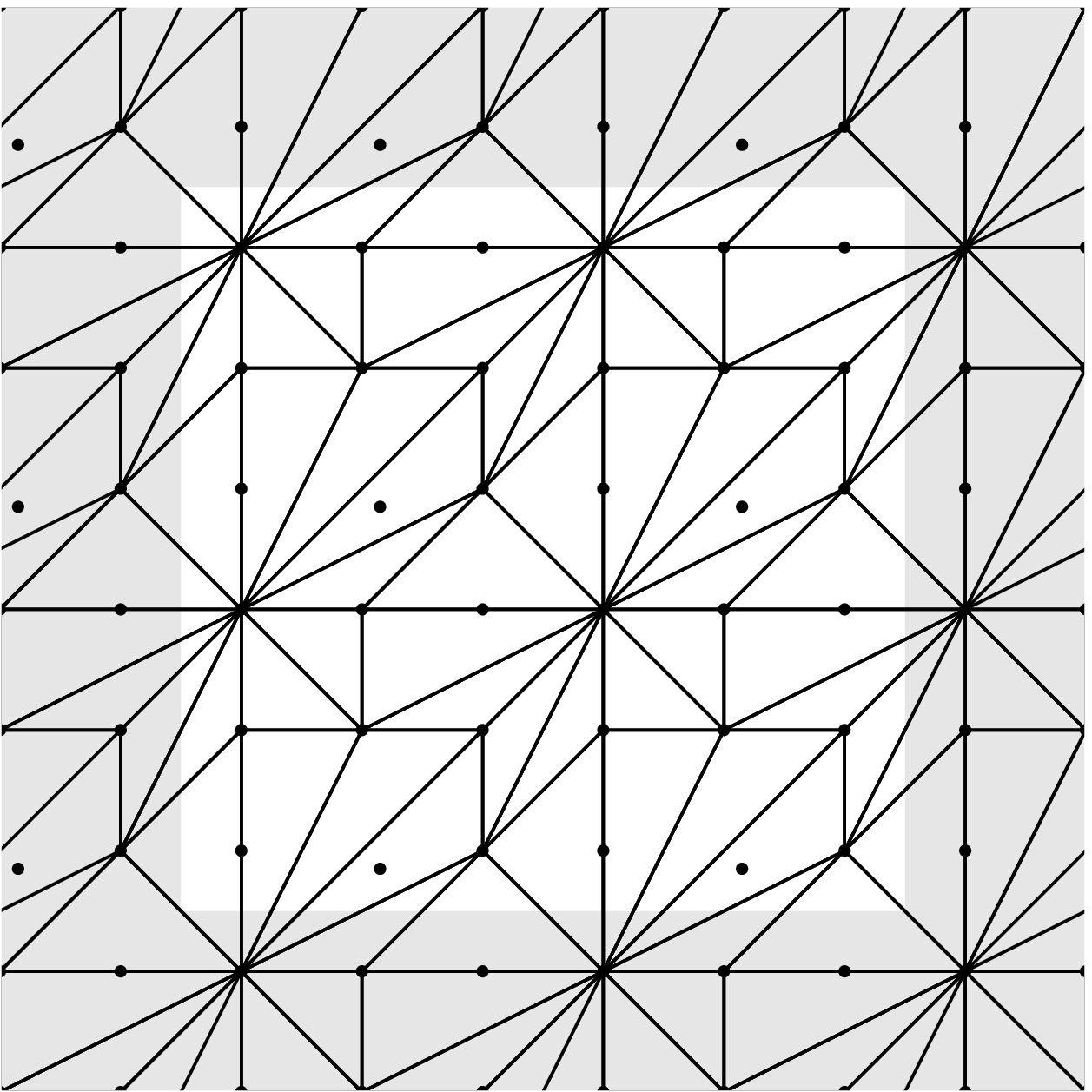}}\\
 (4,8^2) & (3,4,6,4) & (4,6,12) \\ \\
\end{array}\\
\end{array}
\]
\caption{The bond duals of 9 of the Archimedean lattices.}\label{fig_sqlatd}
\end{figure*}

Figure~\ref{fig_sqlat} shows the representations of the
Archimedean lattices $\La$ that we used; the representations
of the planar duals $\Lad$ are shown in Figure~\ref{fig_sqlatd}. The representations
of the site duals $\Lax$ that we used are based on those of the original
lattices $\La$, with extra edges. The planar duals $\Lad$ of
lattices for which the critical probability for bond percolation
is known exactly are omitted.

Let us make some remarks about specific lattices. We wanted our lattices to
have an axis of symmetry at 45 degrees to the horizontal through a
corner of each fundamental region $S_u$: in all the pictures except $(3,12^2)$ this is the
line down and to the right. In the bond dual of the $(4,6,12)$ lattice the
picture looks asymmetric but that is only due to our squashing to make
it fit the square lattice. In other words if we reflect the graph
about a line at 45 degrees to the horizontal we get an isomorphic
graph.

One lattice, $(3^4,6)$, does not have an axis of symmetry, so we ran the
program horizontally and vertically on this lattice. Both
representations are shown and it is easy to see that one is the
reflection of the other about a line 45 degrees to the horizontal.

For efficiency we tried to avoid having holes (vertices of the lattice
not involved in the graph) in our representations: it was not
practical to avoid this for the Kagom\'e and $(3,12^2)$ lattices (and their
site duals) and the bond dual of $(4,6,12)$.

\eject

\bibliography{confint-pre}

\end{document}